\documentclass[11pt]{article} 
\usepackage[utf8]{inputenc}
\RequirePackage{amsthm,amsmath,amsfonts,amssymb,mathtools}
\RequirePackage{graphicx}
\RequirePackage{natbib}

\usepackage{authblk} 
\usepackage[dvipsnames]{xcolor}
\usepackage[plainpages=false,pdfpagelabels,
colorlinks=true,linkcolor=RedOrange,
urlcolor=MidnightBlue,citecolor=MidnightBlue]{hyperref}
\providecommand{\keywords}[1]{\small \textbf{\textit{Keywords---}} #1}
\usepackage{kpfonts}
\usepackage{fullpage}

\usepackage{my-macros}
\usepackage{multirow}
\usepackage{booktabs}

\usepackage{subfig}
\usepackage{tikz}
\usepackage{pgfplots}
\pgfplotsset{compat=1.18} 

\title{Detecting Weak Distribution Shifts via Displacement Interpolation}
\begin{document}

\author[1]{YoonHaeng Hur}
\author[1]{Tengyuan Liang\thanks{\tt Email:\href{mailto:tengyuan.liang@chicagobooth.edu}{tengyuan.liang@chicagobooth.edu}.
Liang acknowledges the generous support from the NSF Career Award (DMS-2042473), and the William Ladany Faculty Fellowship at the University of Chicago Booth School of Business.}}

\affil[1]{University of Chicago}


\maketitle

\begin{abstract}
    Detecting weak, systematic distribution shifts and quantitatively modeling individual, heterogeneous responses to policies or incentives have found increasing empirical applications in social and economic sciences. Given two probability distributions $P$ (null) and $Q$ (alternative), we study the problem of detecting weak distribution shift deviating from the null $P$ toward the alternative $Q$, where the level of deviation vanishes as a function of $n$, the sample size. We propose a model for weak distribution shifts via displacement interpolation between $P$ and $Q$, drawing from the optimal transport theory. We study a hypothesis testing procedure based on the Wasserstein distance, derive sharp conditions under which detection is possible, and provide the exact characterization of the asymptotic Type I and Type II errors at the detection boundary using empirical processes. We demonstrate how the proposed testing procedure works in modeling and detecting weak distribution shifts in real data sets using two empirical examples: distribution shifts in consumer spending after COVID-19, and heterogeneity in the published p-values of statistical tests in journals across different disciplines.
 \end{abstract}


\tableofcontents



\section{Introduction}
\label{sec:intro}

Classic detection problem aims to distinguish a shifted distribution $Q$ from the null distribution $P$ based on data, formulated as a nonparametric goodness-of-fit test:
\begin{equation}
    \label{eq:detection}
        H_0 : X_1, \ldots, X_n \overset{\mathrm{i.i.d.}}{\sim} P 
        \quad \text{vs.} \quad       
        H_1 : X_1, \ldots, X_n \overset{\mathrm{i.i.d.}}{\sim} Q.        
\end{equation}
It is understood that the power of certain test statistics, such as Kolmogorov-Smirnov \citep{kolmogorov1933sulla,smirnoff1939ecarts} or Anderson-Darling \citep{10.1214/aoms/1177729437}, is asymptotically $1$, implying that detection is possible if the sample size is large. 

In many applications, one confronts situations where the signal in the alternative distribution is weaker. One natural formulation is to replace $Q$ in $H_1$ \eqref{eq:detection} by a suitable interpolation scheme between $P$ and $Q$. A common choice is the linear interpolation $(1 - \epsilon) P + \epsilon Q$, or the so-called Huber's $\epsilon$-contamination model \citep{10.1214/aoms/1177703732}, where one parametrizes $\epsilon = \epsilon_n \to 0$ as $n \to \infty$ to represent weak signals. This model represents that only a small $\epsilon$-fraction of the data deviates from $P$, often seen in applications such as large-scale inference with high-throughput measurements \citep{efron_2010}. For instance, microarray data---widely used in genetics and genomics---measure the expression level of thousands of genes, where only some portions are relevant to detecting certain diseases. See \cite{donoho2015higher} for a comprehensive overview of such applications. In those applications, weak signals are modeled as small perturbations in the frequencies of the histogram of the data, that is, deviations along the $y$-axis.

Detecting the presence of distribution shifts is also essential in social and economic applications. There, the main interest is often to quantify individual, heterogeneous responses to policies or other incentives. For instance, consider the economic impacts of the coronavirus pandemic (COVID-19) and government policies \citep{chetty2020economic}; when analyzing how the economy---measured by weekly statistics such as consumer spending---recovers from the shock caused by the pandemic, it is natural to consider individual, heterogeneous shifts to reflect that each individual adjusts the spending differently according to the income level and other characteristics. In this context, the individual responses can be modeled as perturbations along the $x$-axis of the data histogram, whereas the aforementioned perturbation along the $y$-axis (frequencies) rooted in engineering applications is arguably less informative. Another example is the study of the effect of financial or non-financial incentives given to students or parents to improve educational performance, measured by test scores \citep{fryer2015parental,levitt2016behavioralist}. There, the main interest is measuring the shifts in test score distributions in response to the incentives. It is worth noting that the shifts are arguably heterogeneous, depending on gender, race/ethnicity, and other demographic information \citep{levitt2016behavioralist}.

In the above examples, perturbations along the $y$-axis---the fraction of individuals deemed responsive to a given policy or incentive---may be of interest but are insufficient to model individual, heterogeneous responses. Unlike the earlier genetics example, where the identification of a handful of non-null genes showing distinct expression levels provides crucial scientific evidence of certain diseases, the goal of socioeconomic research is beyond simply finding a proportion of shifted observations; the fundamental interest is to quantify individual, heterogeneous responses to policy, which is more relevant to the shift along the $x$-axis of histograms.

This paper proposes a different interpolation scheme for detecting weak distribution shifts motivated by the above. We study a natural test statistic motivated by optimal transport and conduct an exact study of the asymptotic power of the test. We represent the signal strength as the level of deviation, rather than the proportion of deviated units. In a nutshell, we are interested in how much the data histogram shifts along the $x$-axis instead of the $y$-axis. To this end, we model weak signals using displacement interpolation \citep{mccann1997convexity} with a viewpoint from optimal transport theory \citep{Villani_2003}, where the interpolation parameter $\epsilon$ represents a certain transport distance of data from the null $P$ along the geodesics, thus characterizing the level of distribution shifts from $P$.

\subsection{Displacement interpolation and optimal transport}
\label{sec:dis+ot}

Let $P$ and $Q$ be Borel probability measures on $\R$ whose cumulative distribution functions are $F$ and $G$, respectively; also, let $F^{-1}$ and $G^{-1}$ be their quantile functions, respectively. Displacement interpolation between $P$ and $Q$ is defined as
\begin{equation*}
    P_\epsilon := ((1 - \epsilon) \mathrm{Id} + \epsilon T)_\# P \quad \forall \epsilon \in [0, 1],
\end{equation*}
where $\mathrm{Id} \colon \R \to \R$ is the identity map and $T = G^{-1} \circ F$, namely, the composition of $G^{-1}$ and $F$. Here, $S_\# P$ denotes the pushforward measure of $P$ by a map $S \colon \R \to \R$ defined by $S_\# P(B) = P\{x \in \R : S(x) \in B\}$ for any Borel set $B \subset \R$. If $P$ is absolutely continuous with respect to the Lebesgue measure, it is known that $T = G^{-1} \circ F$ is the unique monotone map such that $T_\# P = Q$. More importantly, it is also the unique solution to the following optimal transport problem \citep{Brenier_1991,Villani_2003} provided $P$ and $Q$ have finite second moments, namely,
\begin{equation}
    \label{eq:Monge}
    T = \argmin_{\substack{S \colon \R \to \R \\ S_\# P = Q}} \int_{\R} |x - S(x)|^2 \dd{P(x)}.
\end{equation}
Intuitively, the constraint $S_\# P = Q$ means that a map $S$ transports mass from the source distribution $P$ to the target distribution $Q$ by moving the infinitesimal mass at $x$ in the support of $P$ to $S(x)$ such that $\dd{P(x)} = \dd{Q(S(x))}$. Accordingly, \eqref{eq:Monge} tells that the unique monotone map $T = G^{-1} \circ F$ provides the optimal way of transporting mass from $P$ to $Q$, minimizing the squared transport distance $|x - T(x)|^2$ on average.

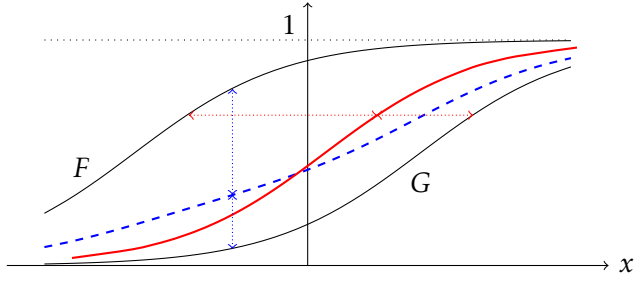
\begin{figure}[ht]
    \centering
    \begin{tikzpicture}[domain=0:4]
        \draw[->] (-4,0) -- (4,0) node[right] {$x$};
		\draw[->] (0,0) -- (0,3.5) ;

        \draw[domain = -3.5:3.5, smooth, variable=\x] plot ({\x},{3/(1+e^(-(\x + 2.3)))});
        \node at (-3, 1.3) {$F$};

		\draw[domain = -3.5:3.5, smooth, variable=\x] plot ({\x},{3/(1+e^(-(\x - 1.5)))});
        \node at (1.5, 1.1) {$G$};

        \draw[domain = -3.5:3.5, smooth, blue, dashed, thick, variable=\x] plot ({\x},{(1/3)*3/(1+e^(-(\x + 2.3))) + (2/3)*3/(1+e^(-(\x - 1.5)))});

        \draw[domain = 0.1:2.9, smooth, red, thick, variable=\u] plot ({(1/3)*(-ln(3/\u-1)-2.3) + (2/3)*(-ln(3/\u-1)+1.5)},{\u});

		\draw[dotted] (-3.5,3) -- (3.5,3);

        \draw[red, densely dotted, <->] (-1.58,2) -- (0.92,2);
        \draw[red, densely dotted, <->] (0.92,2) -- (2.2,2);

        \node[left] at (0,3.2) {1};

        \draw[blue, densely dotted, <->] (-1,0.94) -- (-1,2.35);
        \draw[blue, densely dotted, <->] (-1,0.22) -- (-1,0.94);

    \end{tikzpicture}
	\caption{Illustration of two interpolation schemes. Linear interpolation $(1 - \epsilon) P + \epsilon Q$ vertically combines the cumulative distribution functions $F$ and $G$; namely, its cumulative distribution function (blue, dashed) is $(1 - \epsilon) F + \epsilon G$. Meanwhile, displacement interpolation (red, solid) horizontally combines $F$ and $G$, or equivalently, vertically combines the quantile functions $F^{-1}$ and $G^{-1}$. In other words, the quantile function of $((1 - \epsilon) \mathrm{Id} + \epsilon G^{-1} \circ F)_\# P$ is $(1 - \epsilon) F^{-1} + \epsilon G^{-1}$.}
	\label{fig:dis+mix}
\end{figure}

Now, for each $x$, let us view the segment $\{(1 - \epsilon) x + \epsilon T(x) : \epsilon \in [0, 1]\}$ as the transport path from $x$ to the destination $T(x)$ along its displacement $T(x) - x$. It turns out that transporting mass from the source distribution $P$ along such a path gives rise to a geodesic connecting $P$ and $Q$ in the Wasserstein space \citep{Ambrosio_2005}, namely, $\epsilon = \frac{W_p(P, P_\epsilon)}{W_p(P, Q)}$ for all $\epsilon \in [0, 1]$, where for $p \ge 1$, we denote by $W_p$ the Wasserstein-$p$ distance defined by
\begin{equation*}
    W_p(\mu, \nu) := \left(\int_{0}^{1} |F_{\mu}^{-1}(u) - F_{\nu}^{-1}(u)|^p \dd{u}\right)^{1 / p}
\end{equation*}
for any probability measures $\mu, \nu$ whose quantile functions are $F_{\mu}^{-1}, F_{\nu}^{-1}$, respectively. To see this, it suffices to observe that the quantile function of $P_\epsilon$ is $(1 - \epsilon) F^{-1} + \epsilon G^{-1}$ as $P = (F^{-1})_\# U$ implies $P_\epsilon = ((1 - \epsilon) \mathrm{Id} + \epsilon G^{-1} \circ F)_\# ((F^{-1})_\# U) = ((1 - \epsilon) F^{-1} + \epsilon G^{-1})_\# U$, where $U$ is the Lebesgue measure on $[0, 1]$; see Figure \ref{fig:dis+mix} for details. Therefore, displacement interpolation represents the optimal path---shortest path under the distance $W_p$---for transporting mass from $P$ to $Q$, where the parameter $\epsilon$ naturally characterizes the deviation from $P$ via the relative distance $\epsilon = \frac{W_p(P, P_\epsilon)}{W_p(P, Q)}$.\footnote{It is possible to generalize the aforementioned concepts---displacement interpolation, optimal transport, and the Wasserstein-$p$ distance---to $\R^d$ with $d > 1$; see \citep{Villani_2003}.}

We finish this section with a concrete example to contrast two interpolation schemes. Consider $P = N(0, 1)$ and $Q = N(\tau, \sigma^2)$ for some $\tau, \sigma > 0$. Linear interpolation $(1 - \epsilon) P + \epsilon Q$ is the mixture of two Gaussian distributions $P, Q$, where the parameter $\epsilon$ is essentially the frequency of signals from $Q$. Displacement interpolation is the optimal transport path from $P$ to $Q$, where $T(x) = \sigma x + \tau$ for all $x \in \R$ and $P_\epsilon := ((1 - \epsilon) \mathrm{Id} + \epsilon T)_\# P = N(\epsilon \tau, (1-\epsilon+\epsilon \sigma)^2)$; in words, the optimal path is simply to move each $x$ to $\sigma x + \tau$. The resulting interpolation $N(\epsilon \tau, (1-\epsilon+\epsilon \sigma)^2)$ suggests that the parameter $\epsilon$ provides an intuitive characterization of the signal strength which controls the level of distribution shift. Lastly, it is worth noting that displacement interpolation preserves the unimodality as the interpolation is always Gaussian, whereas linear interpolation may result in two modes.

\subsection{Distribution shift as displacement interpolation}
\label{sec:covid}

As discussed in the previous section, displacement interpolation constructs an optimal path transporting mass from $P$ to $Q$, providing a natural notion of distribution shift from $P$ to $Q$. Such a viewpoint has recently found several applications in machine learning and computer vision, where displacement interpolation serves as a method to optimally synthesize two objects, such as textures \citep{rabin2012wasserstein}, colors \citep{7974883}, and styles \citep{pmlr-v108-mroueh20a}. 

\begin{figure}[!ht]
    \centering
    \subfloat[Daily average consumer spending]{\includegraphics[width=0.9\textwidth]{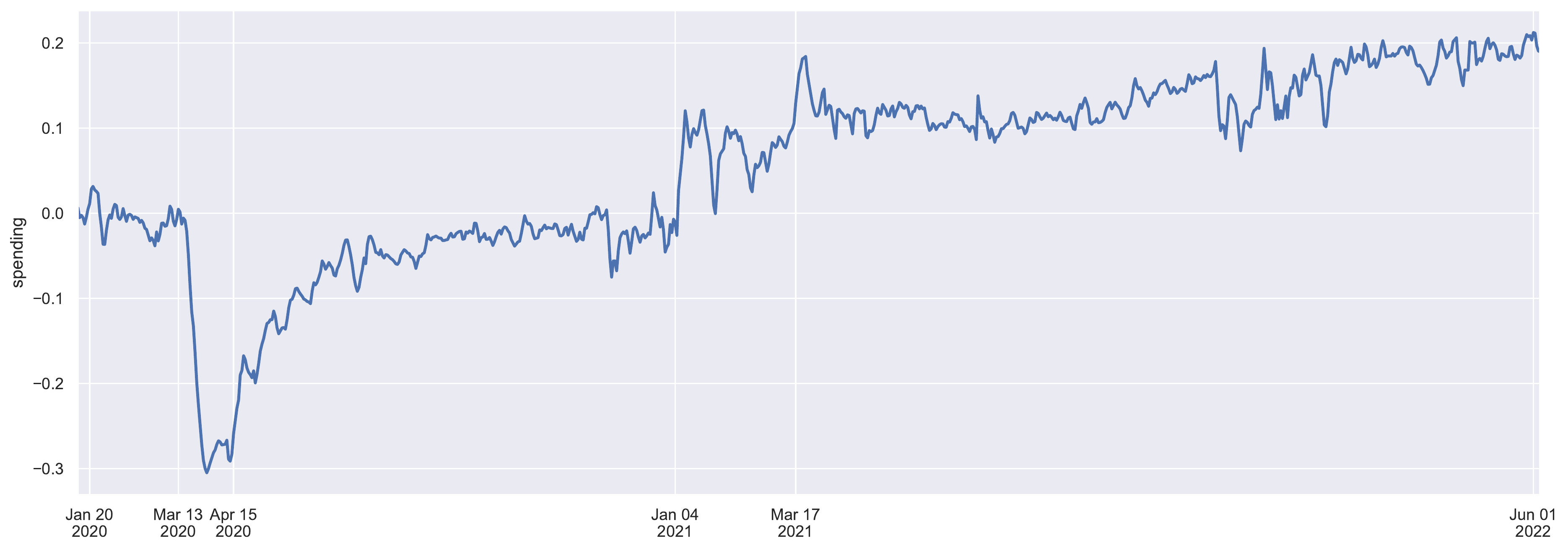}}
    \qquad
    \subfloat[Monthly spending \\ 12 months from March 2020]{\includegraphics[width=0.497\textwidth]{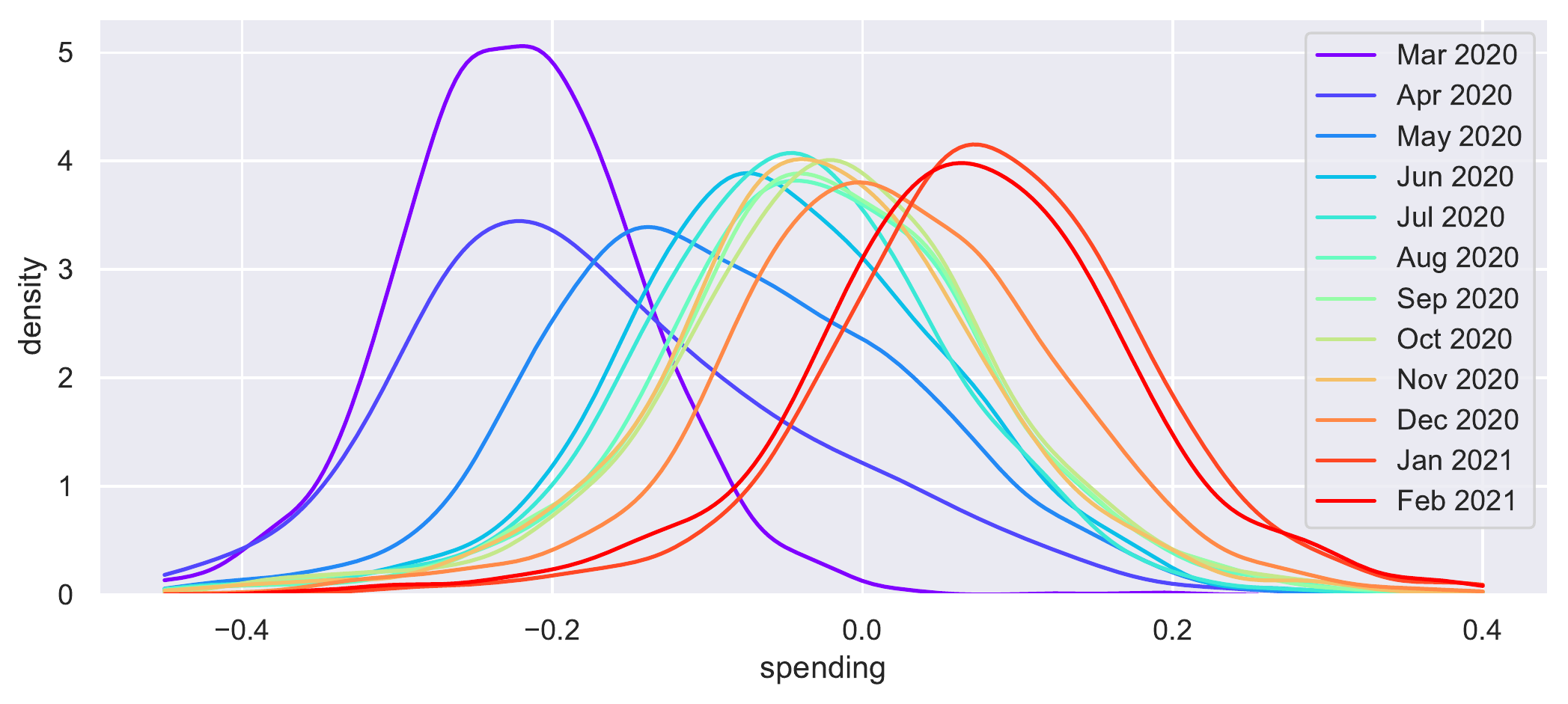}}
    \subfloat[Monthly spending \\ 12 months from March 2021]{\includegraphics[width=0.497\textwidth]{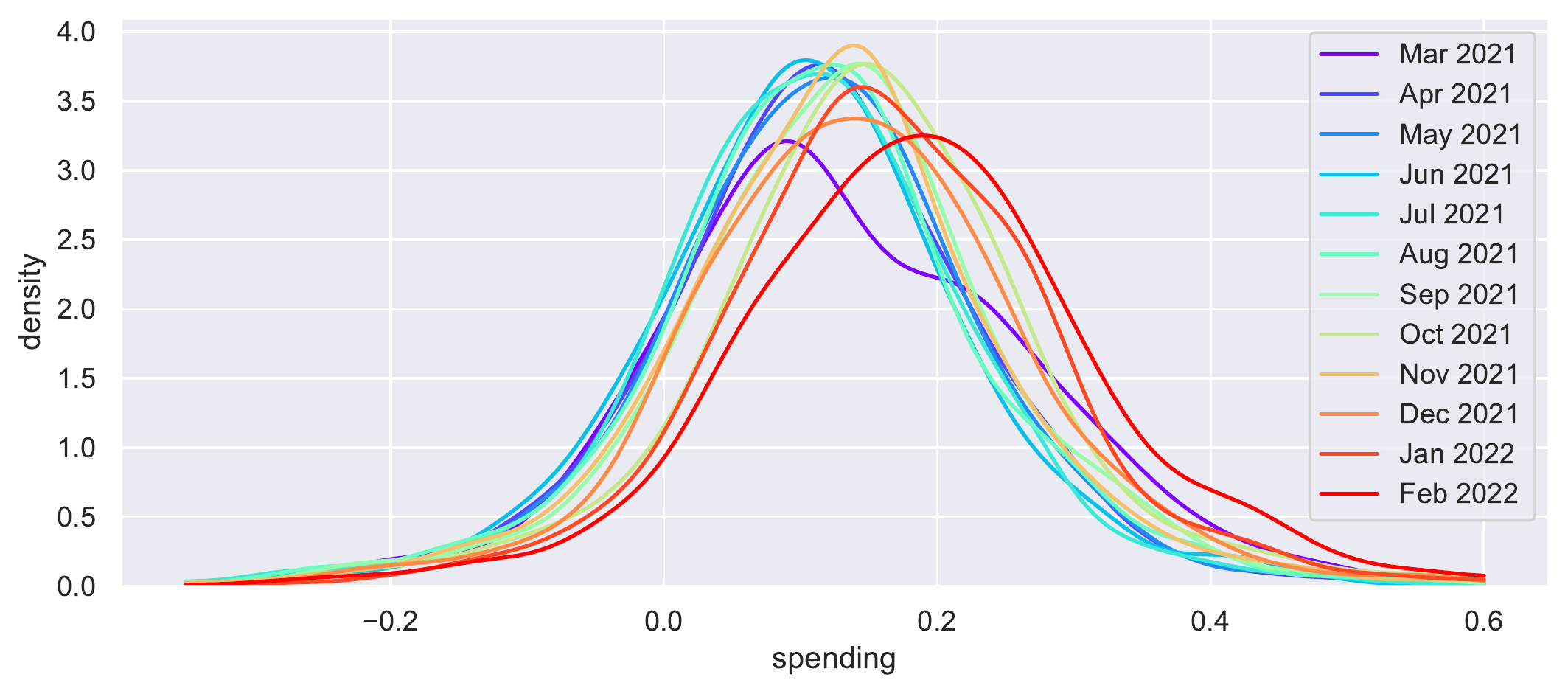}} 
    \caption{\small (a) plots a time series of consumer spending from January 13, 2020 to June 5, 2022, where some noticeable events are marked: March 13, 2020 (national emergency declared) and April 15, 2020, January 4, 2021, and March 17, 2021 (first, second, and third stimulus payments start, respectively). (b) and (c) show the smoothed histograms of monthly average consumer spending by county from March 16, 2020 to March 15, 2021 and from March 16, 2021 to March 15, 2022, respectively.}
    \label{fig:covid_trends}
\end{figure}

Before diving into theory and methodology, in this section, we adopt this viewpoint to model distribution shifts using real-world data; we see that displacement interpolation is better suited to model such distribution shifts than linear interpolation. Later, after laying out the theory, we revisit this empirical example to study the power of our testing procedure. To this end, we utilize the data from \citep{chetty2020economic}, which studies the economic impacts of the coronavirus pandemic (COVID-19) and policy responses in the United States using a wide range of statistics, such as consumer spending, business revenues, employment rates, and so on.\footnote{Data are publicly available at \url{https://tracktherecovery.org/}.} Here, we focus on the consumer spending data---recorded as seasonally adjusted percent changes based on anonymized card transactions data---and analyze how consumer expenditures have recovered from the steep plunge caused by COVID-19. Figure \ref{fig:covid_trends}(a) shows the average consumer spending between January 2021 and June 2022. More specifically, we look at the distribution of monthly average consumer spending over 1,655 counties, where the spending is disaggregated based on the ZIP code where the cardholder lives. The monthly average is calculated by averaging the daily expenditures from the 16th of each month to the 15th of the following month. Figure \ref{fig:covid_trends}(b) shows the smoothed histograms of monthly average spending for the first 12 months since March 2020, which clearly shows how the distribution has shifted in the increasing spending direction, namely, the spending is recovering from the shock of the pandemic. Meanwhile, Figure \ref{fig:covid_trends}(c) shows the next 12 months from March 2021, where we can no longer observe such an evident distribution shift, suggesting the consumer spending has stabilized after the recovery; see also the corresponding period in Figure \ref{fig:covid_trends}(a).

\begin{figure}[!ht]
    \centering
    \subfloat[Relative distances $\epsilon_t$ and $\gamma_t$]{\includegraphics[width=0.497\textwidth]{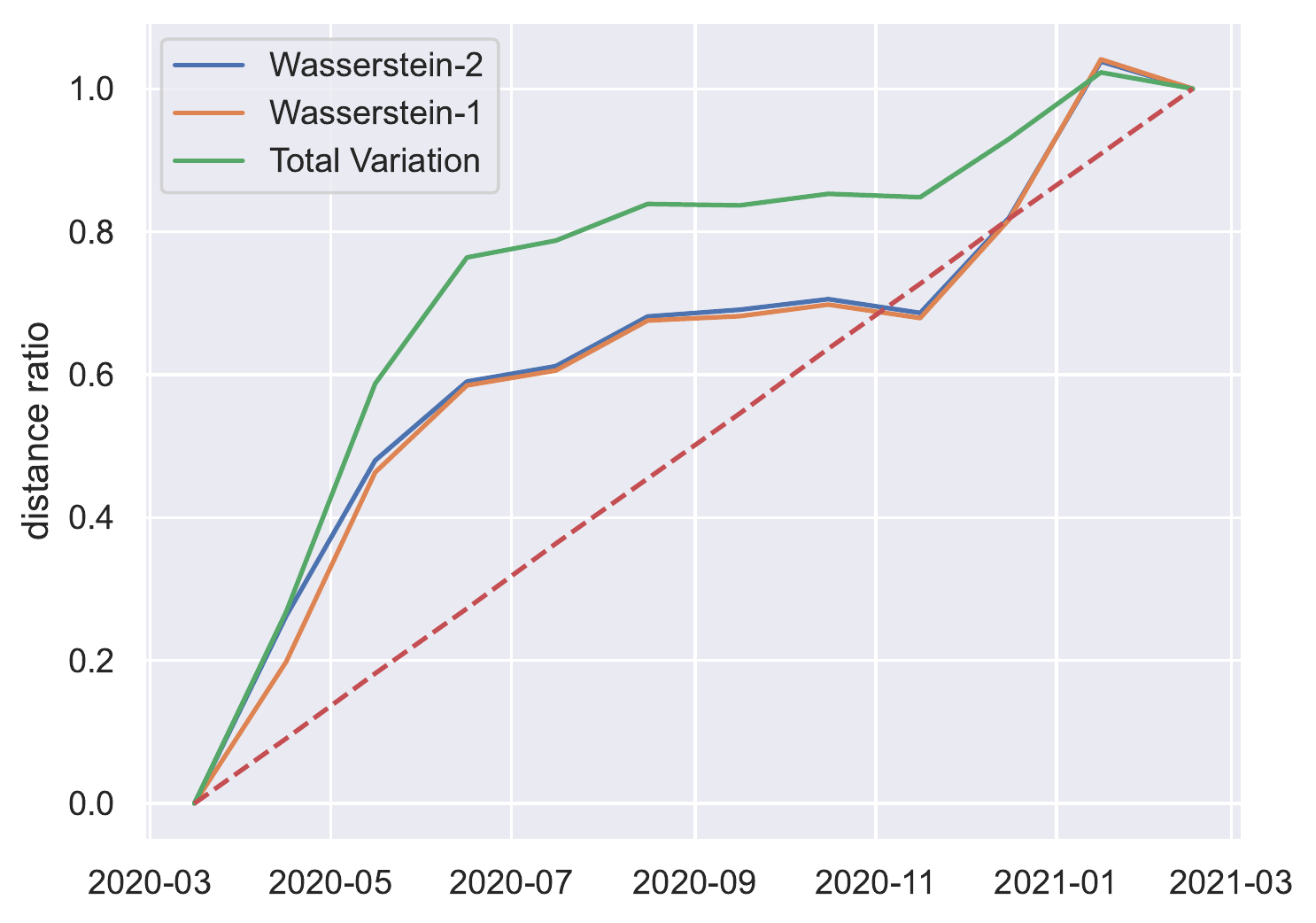}}
    \qquad
    \subfloat[Displacement interpolation $Q_t^{\mathrm{dis}}$]{\includegraphics[width=0.497\textwidth]{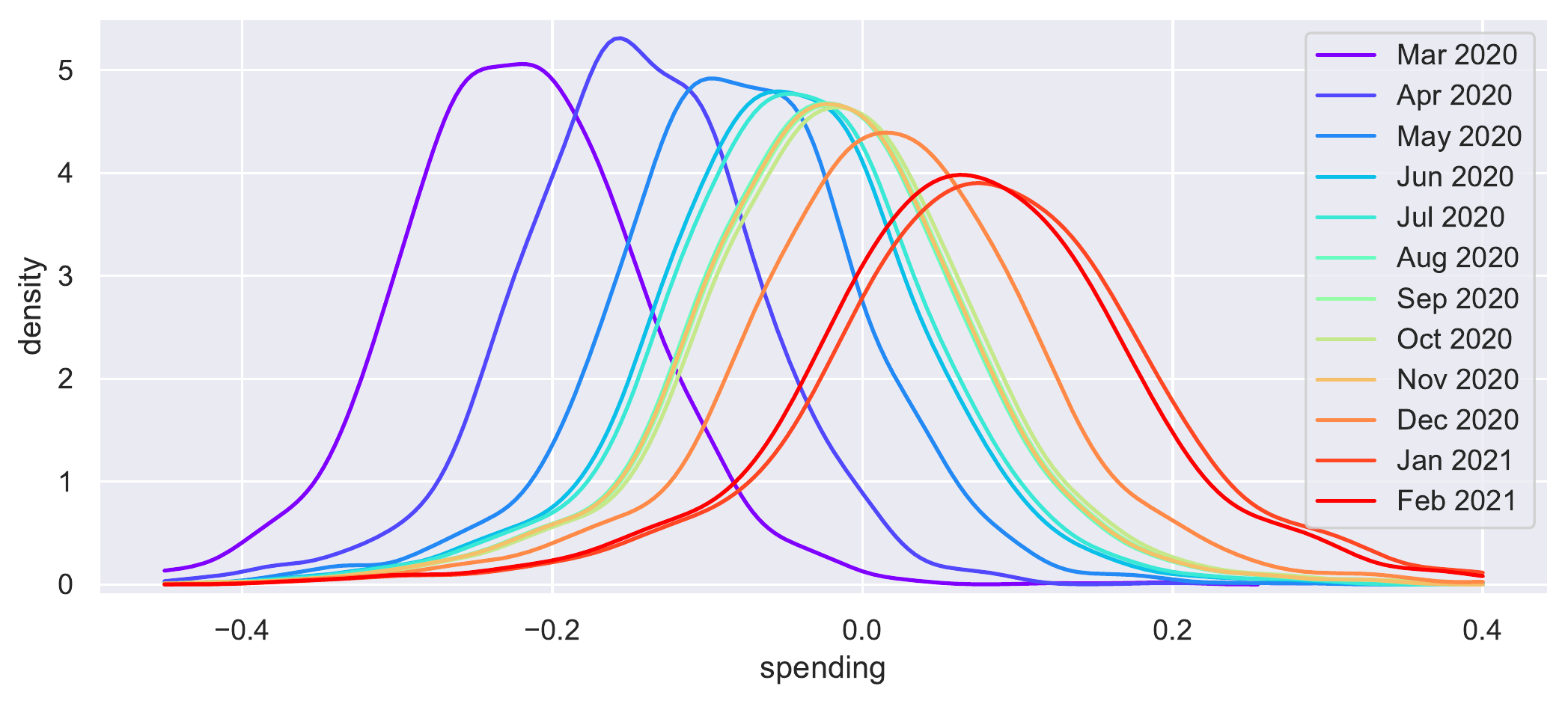}} 
    \subfloat[Linear interpolation $Q_t^{\mathrm{lin}}$]{\includegraphics[width=0.497\textwidth]{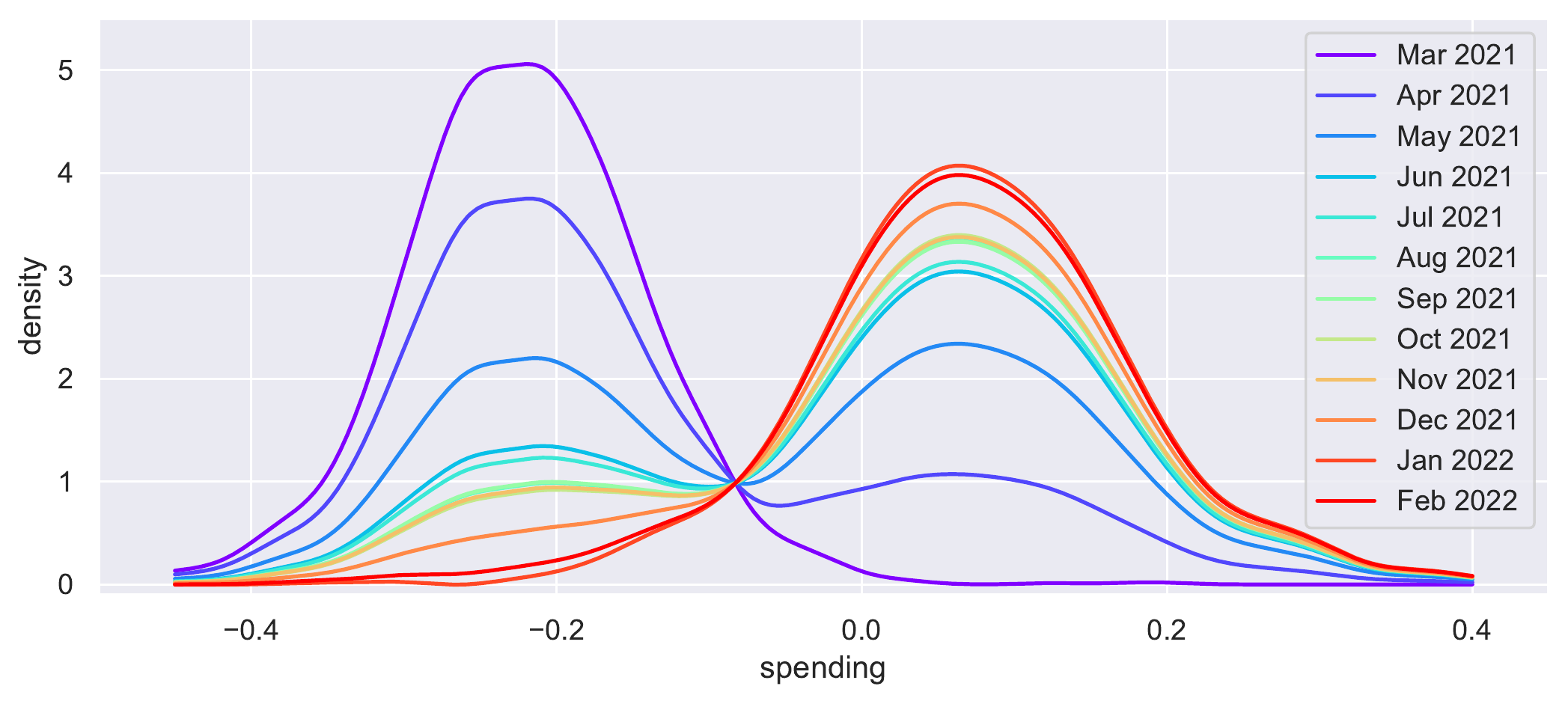}} 
    \caption{\small (a) plots the relative Wasserstein-2 distances $\epsilon_t = \frac{W_2(P_0, P_t)}{W_2(P_0, P_1)}$ and the relative TV distances $\gamma_t = \frac{TV(P_0, P_t)}{TV(P_0, P_1)}$ for $t \in \{0, 1 / 11, \ldots, 1\}$, with the 45 degree line shown as a dashed line; relative Wasserstein-1 distances $\frac{W_1(P_0, P_t)}{W_1(P_0, P_1)}$ are plotted for reference as well, which almost coincide with $\epsilon_t$. (b) and (c) show displacement interpolation $Q_t^{\mathrm{dis}} = ((1 - \epsilon_t) \mathrm{Id} + \epsilon_t T)_\# P_0$ and linear interpolation $Q_t^{\mathrm{lin}} = (1 - \gamma_t) P_0 + \gamma_t P_1$, respectively.}
    \label{fig:covid_distances}
\end{figure}

Empirically, we illustrate that displacement interpolation serves as a reasonable model for the distribution shift during the recovery period shown in Figure \ref{fig:covid_trends}(b). To this end, we generate two interpolation paths, displacement interpolation and linear interpolation, and contrast them with the real data. First, let $\{P_{i / 11}\}_{i = 0}^{11}$ be the distributions of monthly spending by county during that period, namely, they correspond to the 12 histograms of Figure \ref{fig:covid_trends}(b); here, $P_0$ and $P_1$ are the start and end of that period (from March 16, 2020 to April 15, 2020 and from February 16, 2021 to March 15, 2021, respectively). Then, we compute the relative Wasserstein-2 distances $\epsilon_t = \frac{W_2(P_0, P_t)}{W_2(P_0, P_1)}$ and the relative Total Variation (TV) distances $\gamma_t = \frac{TV(P_0, P_t)}{TV(P_0, P_1)}$ for $t \in \{0, 1 / 11, \ldots, 10 /11, 1\}$, as visualized in Figure \ref{fig:covid_distances}(a). From these, we generate displacement interpolation $Q_t^{\mathrm{dis}} = ((1 - \epsilon_t) \mathrm{Id} + \epsilon_t T)_\# P_0$, where $T$ is the composition of the quantile function of $P_1$ and the cumulative distribution function of $P_0$; similarly, we generate linear interpolation $Q_t^{\mathrm{lin}} = (1 - \gamma_t) P_0 + \gamma_t P_1$. Essentially, $Q_t^{\mathrm{dis}}$ amounts to displacement interpolation between $P_0$ and $P_1$ that shifts from $P_0$ at the same rate as $P_t$ under the Wasserstein-2 distance, namely, $W_2(P_0, Q_t^{\mathrm{dis}}) = W_2(P_0, P_t)$; analogously, $Q_t^{\mathrm{lin}}$ corresponds to linear interpolation such that $TV(P_0, Q_t^{\mathrm{lin}}) = TV(P_0, P_t)$. Figure \ref{fig:covid_distances}(b) and Figure \ref{fig:covid_distances}(c) show $Q_t^{\mathrm{dis}}$ and $Q_t^{\mathrm{lin}}$, respectively. Comparing Figure \ref{fig:covid_distances}(b) and Figure \ref{fig:covid_distances}(c) with the real distribution shift in Figure \ref{fig:covid_trends}(b), we can see that displacement interpolation provides a better approximation. Particularly, displacement interpolation preserves the unimodality of distributions as in Figure \ref{fig:covid_trends}(b), while linear interpolation creates two modes, in reminiscence of the Gaussian example in the previous section.

Together with theoretical foundations in Section \ref{sec:dis+ot}, the above example provides the empirical underpinnings of the displacement interpolation model for distribution shifts. Returning to the detection problem \eqref{eq:detection}, in what follows, we will consider a detection problem where $Q$ in $H_1$ is replaced by displacement interpolation $((1 - \epsilon) \mathrm{Id} + \epsilon T)_\# P$ and propose a testing procedure based on the Wasserstein-2 distance. Later, we will apply the proposed testing procedure to revisit the above empirical example and analyze the power under the strong distribution shift during the recovery period.

\subsection{Problem description}
Motivated by the discussion above, we study the problem of detecting weak distribution shifts represented as displacement interpolation: given two distributions $P$ and $Q$ on $\R$ whose cumulative distribution functions are $F$ and $G$, respectively, 
\begin{equation}
    \label{eq:displacement}
    \begin{aligned}
    H_0 : X_1, \ldots, X_n \overset{\mathrm{i.i.d.}}{\sim} P
    \quad \text{vs.} \quad 
    H_1 : X_1, \ldots, X_n \overset{\mathrm{i.i.d.}}{\sim} ((1 - \epsilon) \mathrm{Id} + \epsilon G^{-1} \circ F)_\# P.
    \end{aligned}
\end{equation}
Here, $F^{-1}$ and $G^{-1}$ are the quantile functions of $P$ and $Q$, respectively.

We propose a testing procedure based on the weighted Wasserstein distance between the empirical measure $P_n$---constructed by the observations $X_1, \ldots, X_n$---and the null distribution $P$. Assuming a weak signal in a sense $\epsilon = \epsilon_n \to 0$  as $n \to \infty$, we derive sharp conditions under which detection is possible: (1) when $n^{1/2} \epsilon_n \to 0$, detection is impossible; (2) when $n^{1/2} \epsilon_n \to \infty$, the testing procedure has asymptotic power $1$; (3) at the detection boundary $n^{1/2} \epsilon_n \to \text{constant} \in (0, \infty)$, sharp asymptotic Type I and Type II errors are analyzed using Gaussian processes.

\subsection{Related literature}
\paragraph{Sparse mixture detection} A popular approach to formulating the weak signal detection problem is to replace the alternative hypothesis of \eqref{eq:detection} with Huber's $\epsilon$-contamination model $(1 - \epsilon) P + \epsilon Q$, also known as sparse mixtures, where $\epsilon = \epsilon_n \to 0$ as $n \to \infty$. \cite{donoho2004higher} proposes Tukey's higher criticism as a test statistic and analyzes the asymptotic phase transition depending on the rate $\epsilon = \epsilon_n \to 0$ and the signal strength. A recent endeavor extending the higher criticism test to compare two large frequency tables is given in \cite{donoho2022higher}. See also \cite{tony2011optimal,cai2014optimal} on optimal detection of general sparse mixtures.

\paragraph{Wasserstein distances for testing} 
Our testing procedure is based on the Wasserstein distance, which has been studied extensively in the testing literature. For example, \cite{munk1998nonparametric} and \cite{del2005asymptotics} study Goodness-of-Fit (GoF) testing using the Wasserstein distance between the null hypothesis and the empirical measure based on the observations; see \cite{10.1214/21-EJS1816} for an extension to the multivariate case. As mentioned earlier, the standard GoF testing is related to the classic detection problem \eqref{eq:detection}, where the alternative distribution is fixed to some signal $Q$ or a collection of distributions from a specific parametric family \citep{de2002goodnes,csorgHo2003weighted}.

\paragraph{Notation} For $a, b \in \R$, denote $a \wedge b = \min\{a, b\}$. For positive sequences $\{a_n\}$ and $\{b_n\}$, denote $a_n = o(b_n)$ if $a_n / b_n \to 0$ as $n \to \infty$. We write $Y_n \leadsto Y$ if a sequence $\{Y_n\}$ of random variables converges weakly to some random variable $Y$. 
  
\section{Main Results}

\paragraph{Testing procedure} We propose a distance-based test statistic for the testing problem \eqref{eq:displacement}, motivated by optimal transport and displacement interpolation. More specifically, we compare $P_n$---the empirical measure based on the observations $X_1, \ldots, X_n$---with $P$ under the following distance, called the weighted Wasserstein distance, and reject $H_0$ if it is larger than a specific critical value.

\begin{definition}
    \label{def:weighted_wasserstein}
    Let $\omega$ be a finite Borel measure on $(0, 1)$. For two distributions $\mu$ and $\nu$ on $\R$, we define the $\omega$-weighted Wasserstein distance by 
    \begin{equation*}
        W_{2, \omega}(\mu, \nu) := \left(\int_{0}^{1} |F_\mu^{-1}(u) - F_\nu^{-1}(u)|^2 \dd{\omega(u)}\right)^{1 / 2},
    \end{equation*}
    where $F_\mu^{-1}, F_\nu^{-1}$ denote the quantile functions of $\mu, \nu$, respectively.
\end{definition}

\begin{remark}
    \textup{Note in Definition \ref{def:weighted_wasserstein} that $W_{2, \omega} = W_2$ if $\omega$ is the Lebesgue measure. Weighted versions of the Wasserstein distance are introduced in \citep{de2002goodnes,csorgHo2003weighted}. Here, we introduce the weighted base measure as in the Anderson-Darling test because, for certain detection problems, the signal may hide unevenly among quantiles. For example, the signal is contained in the extreme quantiles in higher criticism \citep{donoho2004higher}.}
\end{remark}

\paragraph{Asymptotic phase transition} We analyze the testing error in the asymptotic regime where $\epsilon = \epsilon_n$ vanishes as $n \to \infty$. More precisely, we rewrite \eqref{eq:displacement} as follows specifying dependency on the sample size $n$:
\begin{equation*}
    \begin{aligned}
    H_0^{(n)} & : X_1, \ldots, X_n \overset{\mathrm{i.i.d.}}{\sim} P, \\
    H_1^{(n)} & : X_1, \ldots, X_n \overset{\mathrm{i.i.d.}}{\sim} ((1 - \epsilon_n) \mathrm{Id} + \epsilon_n G^{-1} \circ F)_\# P,
    \end{aligned}
\end{equation*}
where $\lim_{n \to \infty} \epsilon_n = 0$, namely, $\epsilon_n = o(1)$. Our main result shows that the asymptotic testing error behaves differently depending on the vanishing rate of $\epsilon_n$, highlighting a phase transition phenomenon. In particular, we study a testing procedure such that the asymptotic Type I error is a given level $\alpha \in (0, 1)$ by deriving the limit of the test statistic $W_{2, \omega}(P_n, P)$ under the null hypothesis, based on the fundamental limit theorem of the empirical quantile process: assuming $F$ admits a positive density $f$,  
\begin{equation}
    \label{eq:quantile_convergence_rough}
    \sqrt{n} (F_n^{-1}(u) - F^{-1}(u)) \leadsto \frac{\mathbf{B}_u}{f(F^{-1}(u))},
\end{equation}
where $F_n^{-1}$ is the empirical quantile function based on $X_1, \ldots, X_n$ from $H_0^{(n)}$ and $(\mathbf{B}_u)_{u \in [0, 1]}$ is the standard Brownian bridge, namely, a mean-zero Gaussian process whose covariance satisfies $\mathbb{E} [\mathbf{B}_{u_1} \mathbf{B}_{u_2}] = u_1 \wedge u_2 - u_1 u_2$; rigorous asymptotic analysis is provided in the subsequent section. For such a testing procedure, we can characterize the asymptotic Type II error based on the three phases of $\epsilon_n$ determined by 
\begin{equation*}
    \lim_{n \to \infty} n^{1 / 2} \epsilon_n = \begin{cases}
        0, \\
        \infty, \\
        \gamma \in (0, \infty).
    \end{cases}
\end{equation*}
We formally state the main result as follows.

\begin{theorem}
    \label{thm:asymptotic_errors}
    Suppose $F$ has a density $f$ that is continuous and bounded away from $0$ on some compact interval $I_F$ and is $0$ on $\R \backslash I_F$, $G^{-1}$ is bounded on $(0, 1)$, and $G^{-1} \circ F$ is Lipschitz. Let $\omega$ be a finite Borel measure on $(0, 1)$ that is absolutely continuous with respect to the Lebesgue measure such that $W_{2, \omega}(P, Q) \neq 0$. Also, let $\Psi$ be the cumulative distribution function of 
    \begin{equation}
        \label{eq:limit_distribution_null}
        \int_{0}^{1} \left|\frac{\mathbf{B}_u}{f(F^{-1}(u))}\right|^2 \dd{\omega(u)},
    \end{equation}
    where $(\mathbf{B}_u)_{u \in [0, 1]}$ is the standard Brownian bridge. Fix $\alpha \in (0, 1)$ and let $C_\alpha$ be the $(1 - \alpha)$-th quantile of $\Psi$, that is, $\Psi(C_\alpha) = 1 - \alpha$. Consider the following testing procedure:
    \begin{equation*}
        \text{reject} ~ H_0^{(n)} ~ \text{if and only if} ~~ n W_{2, \omega}^2(P_n, P) > C_\alpha.
    \end{equation*}
    Then, the asymptotic Type I error is $\alpha$, and the asymptotic Type II error is as follows.
    \begin{itemize}
        \item[(i)] If $n^{1 / 2} \epsilon_n \to 0$, the asymptotic Type II error is $1 - \alpha$. 
        \item[(ii)] If $n^{1 / 2} \epsilon_n \to \infty$, the asymptotic Type II error is $0$. 
        \item[(iii)] If $n^{1 / 2} \epsilon_n$ is a constant, say $n^{1 / 2} \epsilon_n = \gamma > 0$, the asymptotic Type II error is
        \begin{equation*}
            \Psi_\gamma(C_\alpha - \gamma^2 W_{2, \omega}^2(P, Q)),
        \end{equation*}
        where $\Psi_\gamma$ is the cumulative distribution function of 
        \begin{equation}
            \label{eq:limit_distribution_alt}
            \int_{0}^{1} \left|\frac{\mathbf{B}_u}{f(F^{-1}(u))}\right|^2 \dd{\omega(u)} + 2 \gamma \int_{0}^{1} \frac{\mathbf{B}_u}{f(F^{-1}(u))} \cdot \left(G^{-1}(u) - F^{-1}(u)\right) \dd{\omega(u)}.
        \end{equation}        
    \end{itemize}
\end{theorem}
\begin{remark}
    \label{rmk:main}
    \textup{Theorem \ref{thm:asymptotic_errors} implies that the asymptotic testing error, namely, the sum of the asymptotic Type I and II errors, is $1$ (undetectable) if $n^{1 / 2} \epsilon_n \to 0$ and $\alpha$ (detectable) if $n^{1 / 2} \epsilon_n \to \infty$. The phase transition occurs at the boundary if $n^{1 / 2} \epsilon_n$ is a constant, where the testing error is determined by the constant $\gamma := n^{1 / 2} \epsilon_n$ and $\Delta := W_{2, \omega}(P, Q)$, which denotes the signal strength. The detection boundary (iii) still holds if we replace $n^{1 / 2} \epsilon_n = \gamma$ with $\lim_{n \to \infty} n^{1 / 2} \epsilon_n = \gamma > 0$; see the formal proof of Theorem \ref{thm:asymptotic_errors} which is deferred to Appendix \ref{sec:proofs}. The main ideas of the analysis will be presented in the next section.}
\end{remark}

\section{Asymptotic Analysis}
In this section, we rigorously derive the asymptotic limit of the test statistic $W_{2, \omega}(P_n, P)$, under the null $H_0^{(n)}$ and the alternative $H_1^{(n)}$, respectively.

\subsection{Preliminaries}
Let $\ell^\infty(0, 1)$ denote the set of all bounded functions defined on $(0, 1)$, which is a Banach space under the uniform norm given by 
\begin{equation*}
    \|h\|_\infty = \sup_{u \in (0, 1)} |h(u)| \quad \forall h \in \ell^\infty(0, 1).
\end{equation*}
We first provide the exact statement of \eqref{eq:quantile_convergence_rough} based on weak convergence in $\ell^\infty(0, 1)$, see \cite{van1996weak}; in what follows, for any random element $Z$ in $\ell^\infty(0, 1)$, let $Z(u)$ denote the random variable at coordinate $u \in (0, 1)$. 

\begin{assumption}
    \label{a:bounded_away}
    $F$ has a density $f$ that is continuous and bounded away from $0$ on some compact interval $I_F$ and is $0$ on $\R \backslash I_F$.
\end{assumption}

\begin{lemma}
    \label{lem:quantile_limit}
    Let $F_n^{-1}$ be the empirical quantile function based on $X_1, \ldots, X_n$ that are i.i.d.\ from a cumulative distribution function $F$. Under Assumption \ref{a:bounded_away}, view $(\sqrt{n} (F_n^{-1} - F^{-1}))_{n \in \N}$ as a sequence of random elements in $\ell^\infty(0, 1)$, then it converges weakly to some tight measurable random element $\mathbf{H}$ in $\ell^\infty(0, 1)$, which we denote as
    \begin{equation}
        \label{eq:quantile_weak_convergence}
        \sqrt{n} (F_n^{-1} - F^{-1}) \leadsto \mathbf{H} ~~ \text{in} ~~ \ell^\infty(0, 1).
    \end{equation}
    The limit $\mathbf{H}$ satisfies the following.
    \begin{itemize}
        \item[(i)] $\{\mathbf{H}(u) : u \in (0, 1)\}$ is a mean-zero Gaussian process with covariance function 
        \begin{equation*}
            \mathbb{E} [\mathbf{H}(u_1) \mathbf{H}(u_2)] = \frac{u_1 \wedge u_2 - u_1 u_2}{f(F^{-1}(u_1)) f(F^{-1}(u_2))} \quad \forall u_1, u_2 \in (0, 1).
        \end{equation*}
        \item[(ii)] The sample path $u \mapsto \mathbf{H}(u)$ is continuous.
    \end{itemize}
\end{lemma}
\begin{remark}
    \label{rmk:on_assumption}
    \textup{Lemma 1 is from Lemma 3.9.23 of \citep{van1996weak}. Assumption \ref{a:bounded_away} ensures that $F^{-1}$ is bounded on $(0, 1)$, thereby viewing $\sqrt{n} (F_n^{-1} - F^{-1})$ as a random element in $\ell^\infty(0, 1)$. Though this assumption rules out distributions supported on the whole real line, such as normal distributions, we can still apply Assumption \ref{a:bounded_away} to such distributions by restricting their support to a sufficiently large yet bounded interval. Alternatively, one may modify Lemma \ref{lem:quantile_limit} by considering weak convergence in $\ell^\infty[\delta, 1 - \delta]$ with a suitable $\delta > 0$. Such a modification provides limit theorems of the integration of $|F_n^{-1} - F^{-1}|^2$ on $[\delta, 1 - \delta]$, often called the trimmed Wasserstein distance \citep{munk1998nonparametric}. It is also possible to modify Lemma \ref{lem:quantile_limit} by considering weak convergence in $L^2(0, 1)$ under suitable assumptions on the behavior of $F^{-1}$ near the endpoints $0$ and $1$ \citep{del2005asymptotics}.}
\end{remark}

\begin{remark}
    \textup{In Lemma \ref{lem:quantile_limit}, the limit $\mathbf{H}$ is tight, meaning that for any $\epsilon > 0$, we can find a compact subset $K$ of $\ell^\infty(0, 1)$ such that $\mathbb{P}(\mathbf{H} \in K) \ge 1 - \epsilon$. Though this technicality is not used explicitly in the main analysis, it is required to apply the extended continuous mapping theorem, which we adapt from Theorem 1.11.1 of \cite{van1996weak} and rewrite as Theorem \ref{thm:extended_continuous_mapping}.}
\end{remark}

Next, we consider the following integrated processes: letting $\mathbf{H}_n := \sqrt{n} (F_n^{-1} - F^{-1})$, define
\begin{align}
    A_n & := \int_{0}^{1} |\mathbf{H}_n(u)|^2 \dd{\omega(u)}, \label{eq:A_n} \\
    B_n & := \int_{0}^{1} |(G^{-1} \circ F)(F^{-1}(u) + n^{-1 / 2} \mathbf{H}_n(u)) - F^{-1}(u)|^2 \dd{\omega(u)}, \label{eq:B_n} \\
    C_n & := \int_{0}^{1} \mathbf{H}_n(u) \cdot \left((G^{-1} \circ F)(F^{-1}(u) + n^{-1 / 2} \mathbf{H}_n(u)) - F^{-1}(u)\right) \dd{\omega(u)}. \label{eq:C_n}
\end{align}
Under the assumptions of Lemma \ref{lem:quantile_limit}, the limit of $A_n$ is essentially the limit of the test statistic $n W_{2, \omega}^2(P_n, P)$ under the null. Later, we will use $B_n$ and $C_n$ to derive the limit of the test statistic under the alternatives. The following lemma derives the limits of the above processes, which we prove in Appendix \ref{sec:proofs}.

\begin{lemma}
    \label{lem:integrated_quantile_convergence}
    Let $F_n^{-1}$ be the empirical quantile function based on $X_1, \ldots, X_n$ that are i.i.d.\ from a cumulative distribution function $F$. Let $\omega$ be a finite Borel measure on $(0, 1)$ that is absolutely continuous with respect to the Lebesgue measure. Under Assumption \ref{a:bounded_away} and assuming $G^{-1}$ is bounded on $(0, 1)$, let $\mathbf{H}_n = \sqrt{n} (F_n^{-1} - F^{-1})$ and $\mathbf{H}$ be the random element mentioned in Lemma \ref{lem:quantile_limit}, then 
    \begin{align}
        A_n & \leadsto \int_{0}^{1} |\mathbf{H}(u)|^2 \dd{\omega(u)}, \label{eq:convergence1} \\
        B_n & \leadsto \int_{0}^{1} |G^{-1}(u) - F^{-1}(u)|^2 \dd{\omega(u)}, \label{eq:convergence2} \\
        C_n & \leadsto \int_{0}^{1} \mathbf{H}(u) \cdot (G^{-1}(u) - F^{-1}(u)) \dd{\omega(u)}, \label{eq:convergence3}
    \end{align}
    where $A_n, B_n, C_n$ are as in \eqref{eq:A_n}, \eqref{eq:B_n}, \eqref{eq:C_n}, respectively.
\end{lemma}

\subsection{Asymptotic distributions}
Now, we analyze the limit of the test statistic $W_{2, \omega}(P_n, P)$ under $H_0^{(n)}$. By \eqref{eq:convergence1} of Lemma \ref{lem:integrated_quantile_convergence}, the following holds.

\begin{proposition}[Limit under the null]
    \label{prop:limit_under_H0}
    For each $n \in \N$, let $P_n$ be the empirical measure based on $X_1, \ldots, X_n$ following $H_0^{(n)}$. Let $\omega$ be a finite Borel measure on $(0, 1)$ that is absolutely continuous with respect to the Lebesgue measure. Under Assumption \ref{a:bounded_away},
    \begin{equation*}
        n W_{2, \omega}^2(P_n, P) \leadsto \int_{0}^{1} |\mathbf{H}(u)|^2 \dd{\omega(u)},
    \end{equation*}
    where $\mathbf{H}$ is the random element mentioned in Lemma \ref{lem:quantile_limit}.
\end{proposition}

Next, we analyze the limit of the test statistic under $H_1^{(n)}$.

\begin{theorem}[Limit under the alternatives]
    \label{thm:asymptotics}
    For each $n \in \N$, let $P_n$ be the empirical measure based on $X_1, \ldots, X_n$ following $H_1^{(n)}$. Let $\omega$ be a finite Borel measure on $(0, 1)$ that is absolutely continuous with respect to the Lebesgue measure. Under Assumption \ref{a:bounded_away} and assuming $G^{-1}$ is bounded on $(0, 1)$ and $G^{-1} \circ F$ is Lipschitz, we can characterize the limit of $n W_{2, \omega}^2(P_n, P)$ as follows.
    \begin{itemize}
        \item[(i)] If $n^{1 / 2} \epsilon_n \to 0$, 
        \begin{equation}
            \label{eq:beta>1/2}
            n W_{2, \omega}^2(P_n, P) \leadsto \int_{0}^{1} |\mathbf{H}(u)|^2 \dd{\omega(u)}.
        \end{equation}
        \item[(ii)] If $n^{1 / 2} \epsilon_n \to \infty$, 
        \begin{equation}
            \label{eq:beta<1/2}
            n^{1 / 2} \epsilon_n^{-1} \left(W_{2, \omega}^2(P_n, P) - \epsilon_n^2 W_{2, \omega}^2(P, Q)\right) \leadsto 2 \int_{0}^{1} \mathbf{H}(u) \cdot \left(G^{-1}(u) - F^{-1}(u)\right) \dd{\omega(u)}.
        \end{equation}
        \item[(iii)] If $n^{1 / 2} \epsilon_n$ is a constant, say $n^{1 / 2} \epsilon_n = \gamma > 0$, 
        \begin{equation}            
            \label{eq:beta=1/2}
            \begin{split}
                n W_{2, \omega}^2(P_n, P) & - \gamma^2 W_{2, \omega}^2(P, Q) \\
                & \leadsto \int_{0}^{1} |\mathbf{H}(u)|^2 \dd{\omega(u)} + 2 \gamma \int_{0}^{1} \mathbf{H}(u) \cdot \left(G^{-1}(u) - F^{-1}(u)\right) \dd{\omega(u)}.
            \end{split}
        \end{equation}
    \end{itemize}
\end{theorem}

\begin{remark}
    \textup{One can also show that the limit distribution in \eqref{eq:beta<1/2} is a mean-zero Gaussian distribution with variance 
    \begin{equation*}
        \int_{0}^{1} \int_{0}^{1} \frac{4(u_1 \wedge u_2 - u_1 u_2)}{f(F^{-1}(u_1)) f(F^{-1}(u_2))} \cdot \left(G^{-1}(u_1) - F^{-1}(u_1)\right) \cdot \left(G^{-1}(u_2) - F^{-1}(u_2)\right) \dd{\omega(u_1)} \dd{\omega(u_2)}.
    \end{equation*}}
\end{remark}
\begin{remark}
    \label{rmk:lipschitzness}
    \textup{The Lipschitzness assumption on $G^{-1} \circ F$ is used for \eqref{eq:beta<1/2}, but not for \eqref{eq:beta>1/2} and \eqref{eq:beta=1/2}. Such an assumption is well studied in the literature; see Appendix of \cite{bobkov2019one}.}
\end{remark}

\section{Simulations}

\subsection{Phase transition and power analysis}
This section delivers results using two experiments to numerically verify Theorem \ref{thm:asymptotic_errors}. For both experiments, we fix $P = \mathrm{Unif}[0, 1]$ and let $\omega$ be the Lebesgue measure so that $\Psi$ is the cumulative distribution function of 
\begin{equation*}
    \int_{0}^{1} |\mathbf{B}_u|^2 \dd{u}.
\end{equation*}
Such a distribution is well studied, and the quantile of $\Psi$ is available, see Table 1 of \cite{10.1214/aoms/1177729437}; we choose $\alpha = 0.05$, then $C_\alpha \approx 0.46136$.

\begin{figure}[ht]
    \centering
	\subfloat[Sum of Type I and Type II errors]{\includegraphics[width=0.497\textwidth]{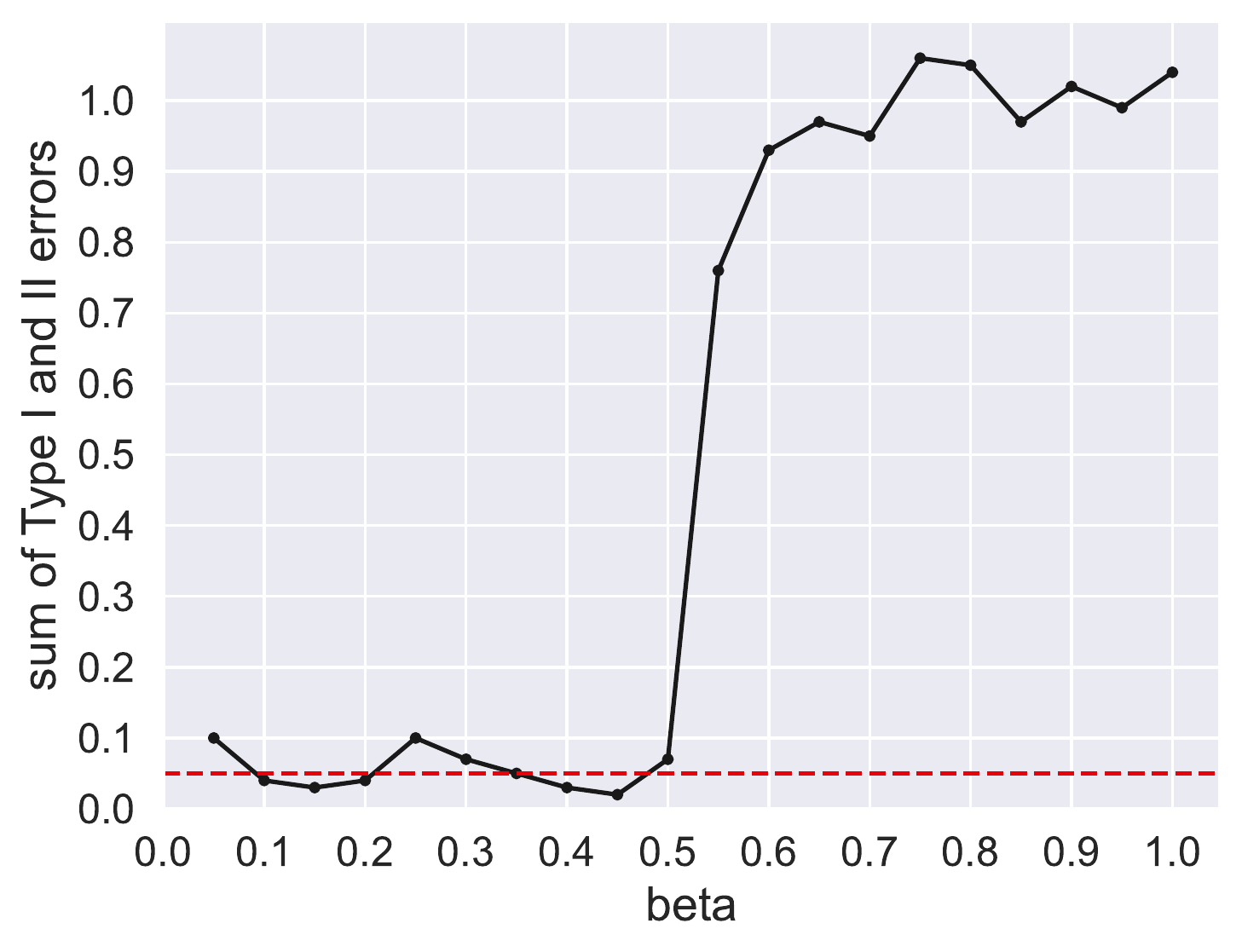}} 
	\subfloat[Type II errors]{\includegraphics[width=0.497\textwidth]{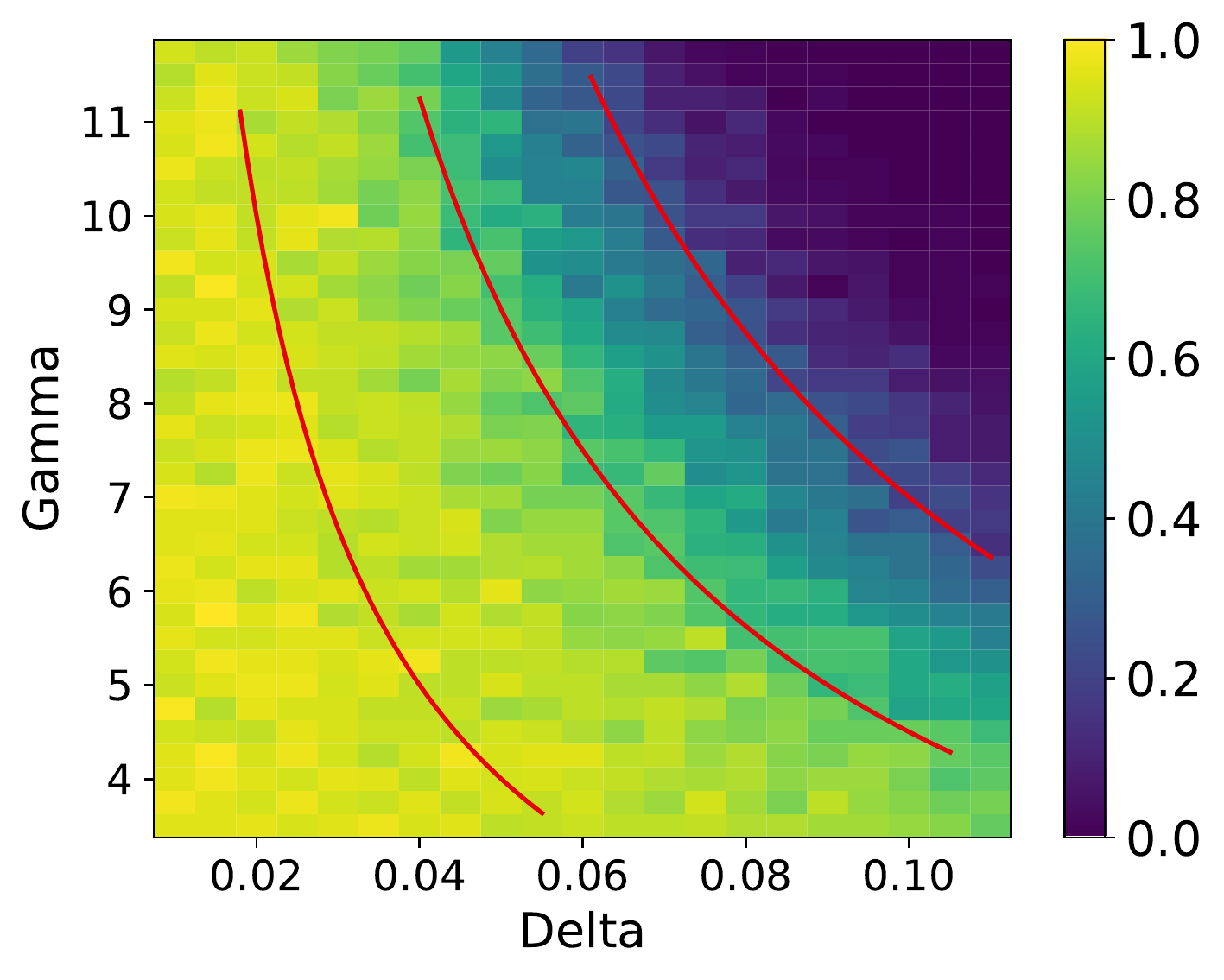}}
	\caption{(a) visualizes the sum of Type I and Type II errors of Experiment 1; the red dashed line represents the level $\alpha = 0.05$. (b) plots Type II errors of Experiment 2 as a color map, where the red solid curves represent $\gamma \cdot \Delta= \text{constant} \in \{0.2, 0.45, 0.7\}$}.
	\label{fig:experiments}
\end{figure}

\paragraph{Experiment 1: phase transition} In this experiment, we verify that the phase transition occurs at $\beta = 0.5$ for the scaling $\epsilon_n = n^{-\beta}$. To this end, we fix $Q = N(0, 1)$ and let $\epsilon_n = n^{-\beta}$ with $n = 10^6$, and a range of $\beta \in \{0.05, 0.1, \ldots, 0.95, 1\}$. For each $\beta$, we compute $n W_2^2(P_n, P)$ using samples from $H_0^{(n)}$ and $H_1^{(n)}$, which we repeat 100 times, then compute the Type I and Type II errors using those 100 realizations, respectively. We plot the sum of Type I and Type II errors against $\beta$. Figure \ref{fig:experiments}(a) confirms that the phase transition occurs at $\beta = 0.5$.

\paragraph{Experiment 2: sharp power analysis} We focus on the power behavior at the detection boundary, namely, $\epsilon_n = \gamma n^{-0.5}$ with $n = 10^6$. Specifically, we analyze the Type II error based on two parameters representing the signal strength: $\gamma$ and $\Delta := W_2(P, Q)$. By Theorem \ref{thm:asymptotic_errors}, the Type II error should be approximately $\Psi_{\gamma}(C_\alpha - \gamma^2 \Delta^2)$. To vary $\Delta$, we parametrized $Q$ as follows: the quantile function of $Q$ is $u \mapsto u + \frac{p}{2 \pi} \sin(2 \pi u)$, where we vary $p \in (0, 1)$; note that this is a valid quantile function as it is monotonically increasing. Then,
\begin{equation*}
    W_2^2(P, Q) = \frac{p^2}{4 \pi^2} \int_0^1 |\sin(2 \pi u)|^2 \, \mathrm{d} u = \frac{p^2}{8 \pi^2} = \Delta^2,
\end{equation*}
which results in $\Delta^2 \in (0, (8 \pi^2)^{-1})$. We can easily generate $Q$'s with desired signal strength $\Delta$ using this parametrization. For each pair $(\Delta, \gamma)$, where $\Delta \in \{0.01, 0.015, \ldots, 0.105, 0.11\}$ and $\gamma \in \{3.5, 3.75, \ldots, 11.5, 11.75\}$, we compute the Type II error and visualize it as a color map. Figure \ref{fig:experiments}(b) visualizes Type II errors on the $\Delta$-$\gamma$ plane using colors; this essentially plots $\Psi_{\gamma}(C_\alpha - \gamma^2 \Delta^2)$. Notice that the level set of Type II errors and the curve $\gamma \cdot \Delta = \text{constant}$ do not perfectly coincide because Type II errors can vary on the curve $\gamma \cdot \Delta = \text{constant}$, namely, $\Psi_{\gamma}(C_\alpha - \text{constant}^2)$ can change as $\gamma$ varies.

\subsection{Comparison to other methods and robustness checks}

This section first compares the power of the proposed procedure and the Kolmogorov-Smirnov (KS) test. Later, we carry out preliminary robustness checks on how the choice of the weight measure $\omega$ affects the power. Throughout this section, we again fix $P = \mathrm{Unif}[0, 1]$ and $n = 10^6$.

\begin{figure}[!ht]
    \centering
    \subfloat[Powers (Example I)]{\includegraphics[width=0.497\textwidth]{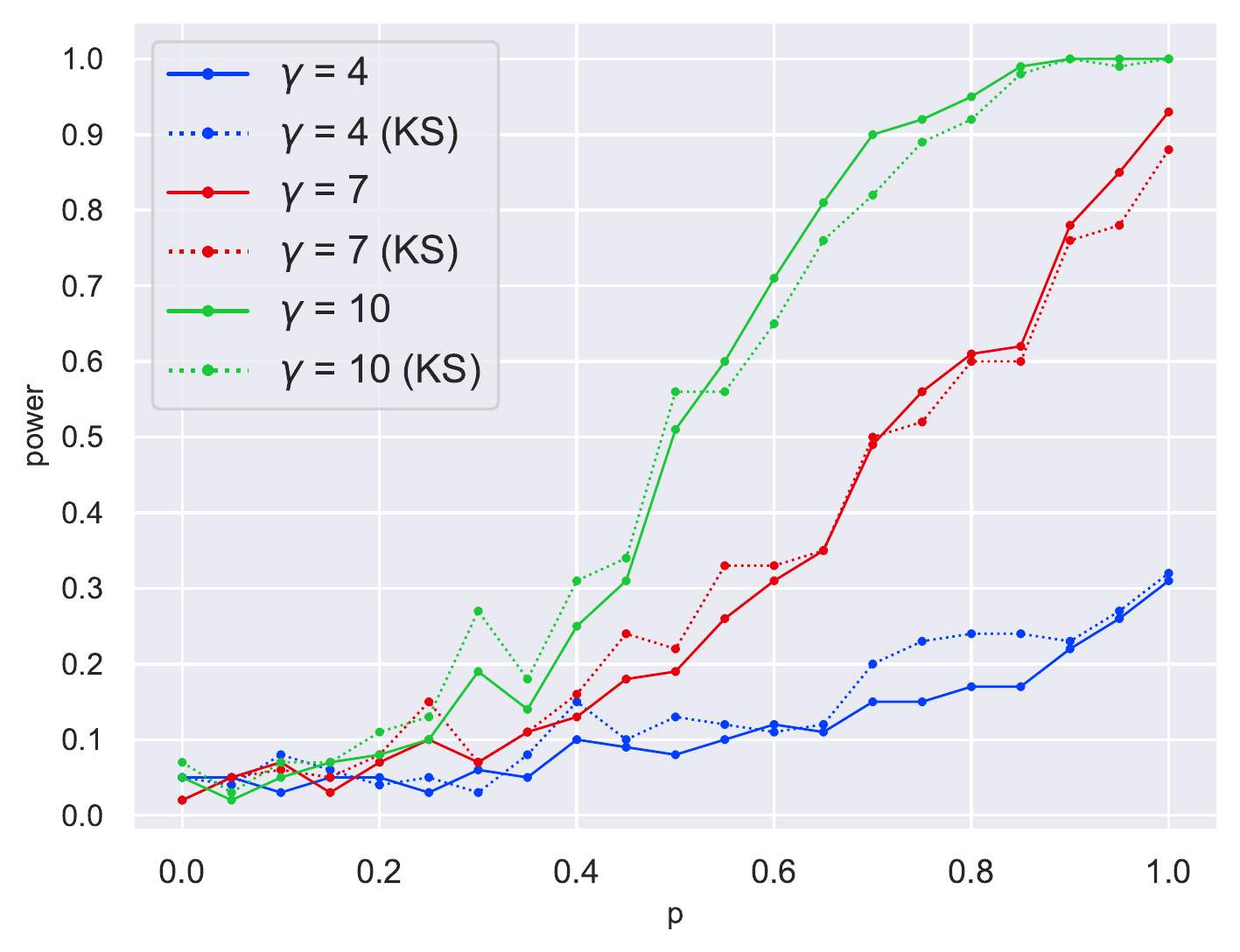}}
    \subfloat[Powers (Example II)]{\includegraphics[width=0.497\textwidth]{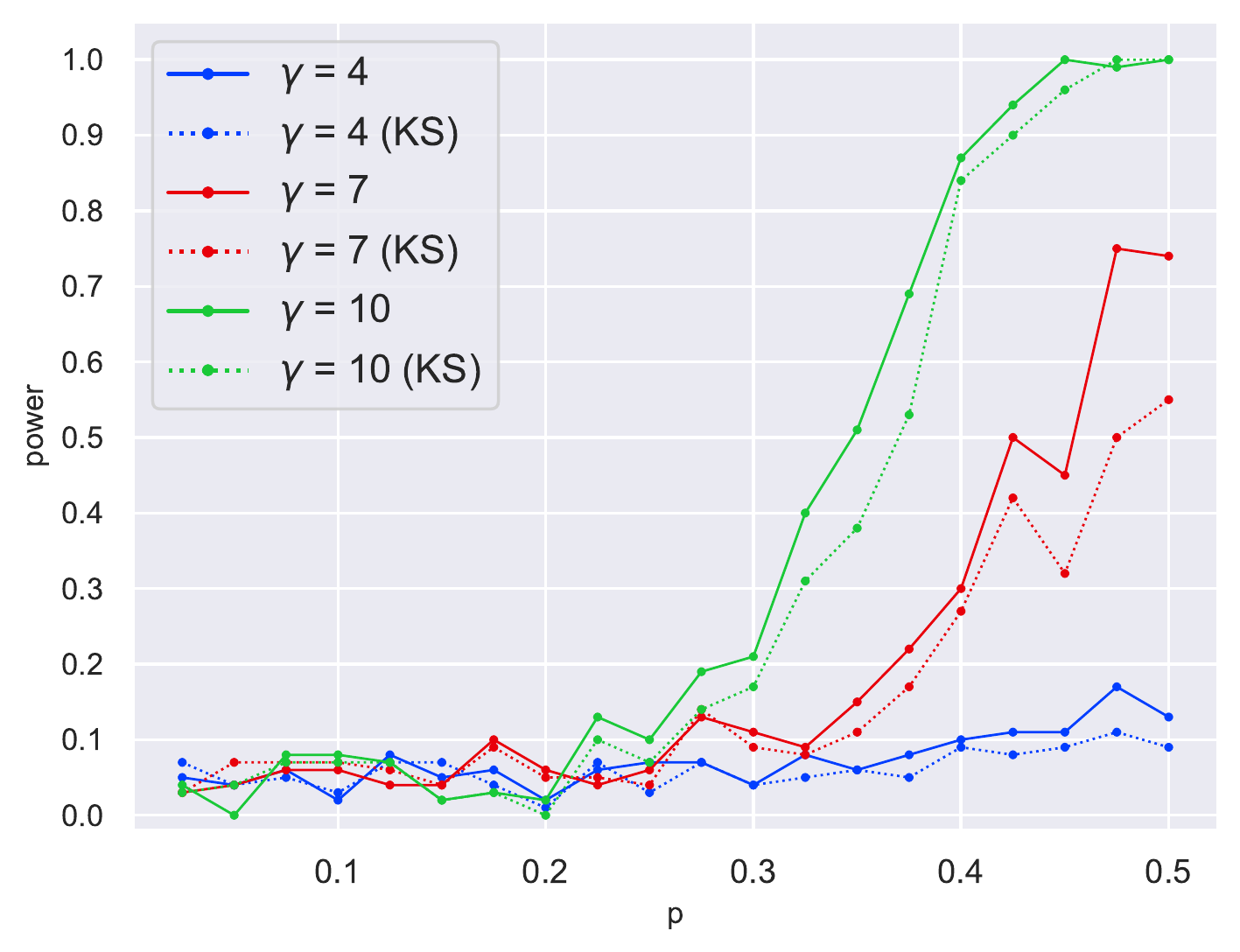}} 
    \qquad
    \subfloat[Powers (Example II) of different weight measures]{\includegraphics[width=0.497\textwidth]{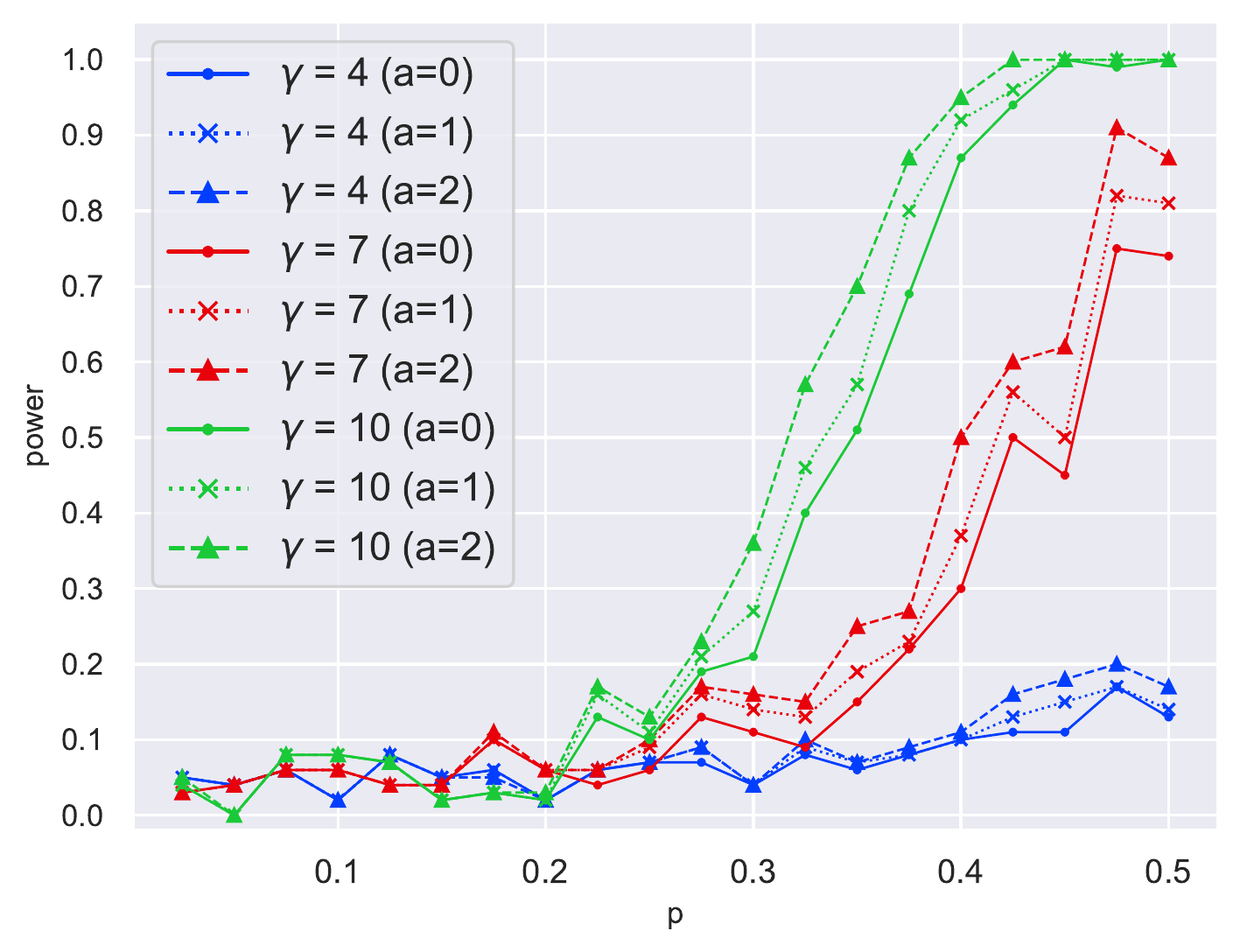}} 
    \caption{\small (a) and (b) show the powers of the proposed testing procedure and the KS test in Example I and Example II, respectively; solid lines correspond to the proposed testing procedure, whereas dotted lines represent the KS test. (c) shows the powers of the proposed testing procedures in the setting of Example II, using the weight measure given by \eqref{eq:quadratic_weight}. By design, the solid lines of (b) and (c) coincide; both correspond to the Lebesgue measure, namely, $a = 0$ in \eqref{eq:quadratic_weight}.}
    \label{fig:comparisons}
\end{figure}

\paragraph{Example I} First, we repeat the previous setting used for the boundary case: $\epsilon_n = \gamma n^{-0.5}$ and the quantile function of $Q$ is $u \mapsto u + \frac{p}{2 \pi} \sin(2 \pi u)$. This time, we parametrize $p \in \{0, 0.05, \ldots, 0.95, 1\}$ instead of $\Delta = W_2(P, Q)$ and restrict our interest to $\gamma \in \{4, 7, 10\}$. For each pair $(p, \gamma)$, we compute the power of the proposed testing procedure---with $\omega$ being the Lebesgue measure---and the KS test; to ensure both are asymptotically level $\alpha$ tests with $\alpha = 0.05$, we use the critical value $C_\alpha \approx 0.46136$ for the proposed testing procedure as before and the critical value $1.36$ for the KS test statistic $\sqrt{n} \cdot \mathrm{KS}(P_n, P) = \sqrt{n} \cdot \sup_{x \in \R} |F_n(x) - F(x)|$ based on Table 1 of \citep{shorack_wellner}. Figure \ref{fig:comparisons}(a) shows the results. First, observe that for $\gamma \in \{7, 10\}$, the powers of both tests exceed $0.5$ once the parameter $p$ is above certain values, which means that both are better than a trivial test that randomly rejects the null with probability $0.5$. In this case, the power of the proposed procedure is slightly larger than that of the KS test, as the solid lines are above the dotted lines. For $\gamma = 4$, the power of the KS test is mostly larger than that of the proposed procedure, as shown in the blue lines; in this case, however, the powers of both tests are far less than $0.5$, implying that both tests practically failed. In other words, the signal strength $\gamma = 4$ is too weak for both tests to detect.

\paragraph{Example II} Next, we consider a different setting where the quantile function of $Q$ differs with that of $P = \mathrm{Unif}[0, 1]$ at tails. To this end, suppose the quantile function of $Q$ is as follows:
\begin{equation}
    \label{eq:quantile_tail_deviated}
    \begin{cases}
        u + 0.45 \cdot \frac{2 p}{\pi} \cos\left(\frac{\pi u}{2 p}\right) & 0 \le u \le p, \\
        u & p \le u \le 1 - p, \\
        u - 0.45 \cdot \frac{2 p}{\pi} \cos\left(\frac{\pi (1 - u)}{2 p}\right) & 1 - p \le u \le 1,
    \end{cases}
\end{equation}
where we parametrize $p \in \{0.025, 0.05, \ldots, 0.5\}$. Note that the quantile function of $Q$ deviates from the identity only at tails, namely, $[0, p] \cup [1 - p, 1]$. As in Example I, we keep using $\epsilon_n = \gamma n^{-0.5}$ for $\gamma \in \{4, 7, 10\}$ and compute the powers for each pair $(p, \gamma)$ which are shown in Figure \ref{fig:comparisons}(b). For $\gamma \in \{7, 10\}$, the results are similar to Example I; for sufficiently large $p$, both tests have enough powers. Now, for $\gamma = 4$, the power of the proposed procedure is slightly larger than that of the KS test; however, as in Example I, both tests are practically useless as their powers are too small; namely, $\gamma = 4$ is again too weak for them to detect.

\paragraph{The role of the weight measure} Lastly, we compare the powers by differing the weight measure $\omega$. Let us keep considering the setting of Example II. As the quantile function $Q$ mainly deviates from that of $P$ at tails, it is reasonable to put more weights at tails for better performance when computing the test statistic. To verify this, consider a simple weight measure whose density with respect to the Lebesgue measure is a quadratic function: for $a \ge 0$,
\begin{equation}
    \label{eq:quadratic_weight}
    \frac{\dd{\omega(u)}}{\dd{u}} = a \left(u - \frac{1}{2}\right)^2 + 1 - \frac{a}{12},
\end{equation}
which satisfies $\int_{0}^{1} \dd{\omega} = 1$. Clearly, $a = 0$ corresponds to the Lebesgue measure; large $a$ means more weights at both tails. We compute the powers of the proposed procedure for $a \in \{0, 1, 2\}$. Here, we again set critical values to ensure that all of them are asymptotically level $\alpha = 0.05$; for $a = 0$, we can reuse the previous critical value as before. For $a \in \{1, 2\}$, we estimate the quantile of the asymptotic distribution $\int_{0}^{1} |\mathbf{B}_u|^2 \dd{\omega(u)}$ by Monte Carlo simulations. Figure \ref{fig:comparisons}(c) shows the results. As expected, we obtain larger powers with the quadratic weight measure, namely, both $a = 1$ and $a = 2$ achieve better performance than the unweighted procedure $a = 0$; particularly, assigning more weights at tails ($a = 2$) yields the largest power. For $\gamma = 4$, though we have increased powers for the weighted procedures, they are still practically too weak to detect the signal, as discussed in the previous examples.

\subsection{Application I: distribution shifts in consumer spending}
We revisit the data example in Section \ref{sec:covid} and apply the proposed testing procedure to study the power. Recall that $\{P_{i / 11}\}_{i = 0}^{11}$ are the distributions of monthly average spending by county during the recovery period between March 16, 2020 and March 15, 2021. For $t \in \{0, 1 / 11, \ldots, 10 / 11, 1\}$ and $n \in \{10, 50, 100, 500\}$, we construct the empirical measure $P_t^n$ using $n$ points that are i.i.d.\ from $P_t$ and compute $W_2(P_0, P_t^n)$. Essentially, for each $t$ and sample size $n$, we want to distinguish $P_t$ from $P_0$ using a finite sample. To this end, we first estimate the quantile of $W_2(P_0, P_0^n)$ via Monte Carlo simulations to define a level $\alpha = 0.05$ test; this will give us a critical value, say, $C_\alpha^n$. Then, for each $t > 0$, we compute $W_2(P_0, P_t^n)$ and reject the null---that is, data are from $P_0$---if it exceeds $C_\alpha^n$; by repeating this for 100 times, we can evaluate the power of the test using simulations. The results are shown in the first 4 rows (Recovery Period) of Table \ref{tab:powers}. In this case, for $t \ge 3 / 11$, namely, after three months from March, 2020, we can detect the distribution shift from $P_0$ for any sample size $n$. In other words, we can tell there has been a significant recovery after three months; indeed, the shift is significant enough so we can tell the difference from $P_0$ by estimating $P_t$ using only 10 randomly chosen counties instead of the total 1,655 counties. 
The first month after $P_0$, that is, for $t = 1 / 11$, the shift is relatively weaker compared to the subsequent months, so $n = 10$ yields power $0.56$, which is not enough to detect the shift; in other words, for $t = 1 / 11$, we need at least $n = 50$ to distinguish it from $P_0$.

How about the power for the stable period? We repeat the above procedure by taking $\{P_{i / 11}\}_{i = 0}^{11}$ as the distributions during the stable period, namely, they amount to the histograms in Figure \ref{fig:covid_trends}(c); again, $P_0$ and $P_1$ correspond to the start and end of this period (from March 16, 2021 to April 15, 2021 and from February 16, 2022 to March 15, 2022). Recall from Figure \ref{fig:covid_trends}(c) that we no longer see a significant distribution shift as in the recovery period. The resulting powers are shown in the last 4 rows (Stable Period) of Table \ref{tab:powers}. Unlike the recovery period, we can observe the powers are much smaller for $n \le 100$; particularly, most of the powers---except for a few periods with $n = 100$---are far below $0.5$, meaning that we cannot detect a meaningful shift from $P_0$ using finite samples. For detection to be possible, we can see that we need a sufficiently large sample size $n = 500$. 

In summary, the results in Table \ref{tab:powers} essentially demonstrate that the proposed testing procedure can detect the distribution shifts during the recovery period as it is reasonably approximated by displacement interpolation as seen in Figure \ref{fig:covid_trends}(b), while there are no such shifts to be captured by the proposed procedure during the stable period because the distributions are similar to each other as shown in Figure \ref{fig:covid_trends}(c).

\renewcommand{\arraystretch}{1.0}
\begin{table}[ht]
	\centering
    \scalebox{0.9}{
	\begin{tabular}{|c|*{11}{c}|r|} 
        \hline
        & \multicolumn{11}{c|}{$11 t$} & \\ 
        \cline{2-12}
		$n$ & 1 & 2 & 3 & 4 & 5 & 6 & 7 & 8 & 9 & 10 & 11 &  \\
		\hline
        10 & 0.56 & 0.98 & 1 & 1 & 1 & 1 & 1 & 1 & 1 & 1 & 1 & \multirow{4}{0.8in}{\textbf{Recovery \\ Period}} \\
        50 & 0.98 & 1 & 1 & 1 & 1 & 1 & 1 & 1 & 1 & 1 & 1 & \\
		100 & 1 & 1 & 1 & 1 & 1 & 1 & 1 & 1 & 1 & 1 & 1 & \\
		500 & 1 & 1 & 1 & 1 & 1 & 1 & 1 & 1 & 1 & 1 & 1 & \\
		\cline{1-13}
		10 & 0.04 & 0.05 & 0.09 & 0.1 & 0.05 & 0.07 & 0.1 & 0.05 & 0.07 & 0.11 & 0.23 & \multirow{4}{0.8in}{\textbf{Stable \\ Period}} \\
        50 & 0.2 & 0.23 & 0.36 & 0.38 & 0.22 & 0.18 & 0.1 & 0.14 & 0.09 & 0.4 & 0.75 & \\
        100 & 0.57 & 0.45 & 0.66 & 0.66 & 0.45 & 0.29 & 0.43 & 0.33 & 0.16 & 0.76 & 0.96 & \\
        500 & 1 & 0.99 & 1 & 1 & 1 & 0.94 & 1 & 0.99 & 0.88 & 1 & 1 & \\
		\hline
	\end{tabular}
    }
	\caption{Powers for the recovery and stable periods.}
    \label{tab:powers} 
\end{table}

\subsection{Application II: p-value heterogeneity across disciplines}
\label{sec:simulation_real}
We apply our testing procedure to the data set from \cite{head2015extent}, which collects p-values of statistical tests published in journals across different disciplines. In this example, we use our procedure to determine if the distribution of p-values differs depending on disciplines.

\begin{figure}[!ht]
    \centering
    \includegraphics[width=0.497\textwidth]{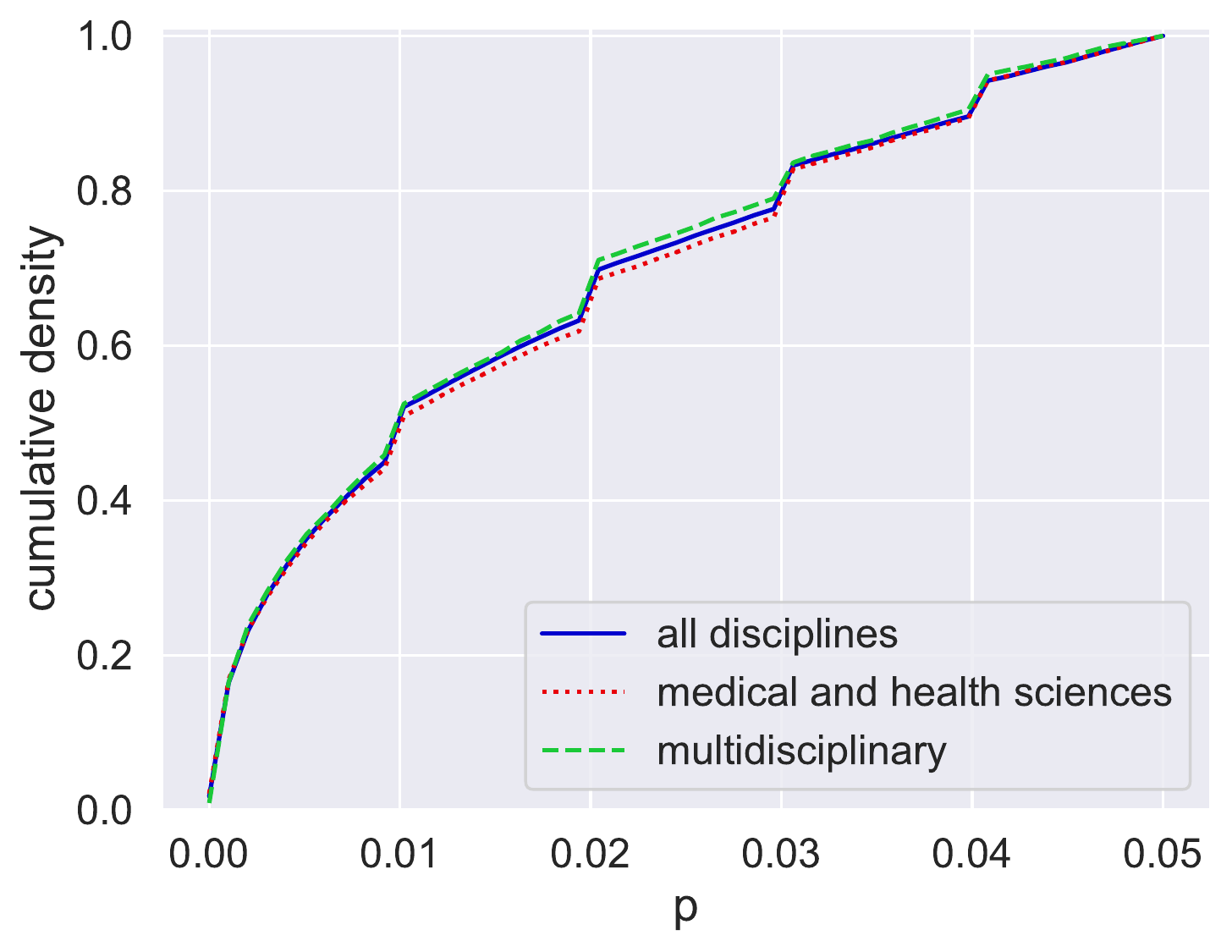}
	\caption{\small This figure plots the cumulative distribution functions of p-values (below $0.05$) across nine disciplines subsampled from \citep{head2015extent}, which is shown as the solid blue curve, along with two disciplines: ``medical and health sciences'' (red, dotted) and ``multidisciplinary'' (green, dashed).}
	\label{fig:p-values}
\end{figure}

To this end, we first collect p-values across 9 disciplines subsampled from total $21$ disciplines, then take $80\%$ of them as the null distribution $P_{\text{Null}}$, which consists of 90,950 observations; this is shown as the solid blue curve in Figure \ref{fig:p-values}. Then, from the remaining $20\%$, we sample two disciplines ``medical and health sciences'' ($P_{\text{MH}}$) and ``multidisciplinary'' ($P_{\text{Mu}}$), which consist of 12,928 and 5,731 observations, respectively.\footnote{They have the largest number of observations among all 9 disciplines.} Figure \ref{fig:p-values} shows the cumulative distribution functions of the null $P_{\text{Null}}$, $P_{\text{MH}}$, and $P_{\text{Mu}}$, which are close to each other. This suggests that both $P_{\text{MH}}$ and $P_{\text{Mu}}$ are weak signals that are hard to distinguish from $P_{\text{Null}}$ visually. To tackle this problem, we apply our testing procedure to $H_0 : P_{\text{MH}} =  P_{\text{Null}}$. First, we compute the test statistic $T_{\text{MH}} := n W_2^2(P_{\text{MH}}, P_{\text{Null}})$, where $n := |P_{\text{MH}}| =$ 12,928 observations. Then, setting $\alpha = 0.05$, we estimate the $(1 - \alpha)$-th quantile of $n W_2^2(P_n, P_{\text{Null}})$, where $P_n$ is the empirical measure based on $X_1, \ldots, X_n$ that are i.i.d.\ from $P_{\text{Null}}$, which can be done by resampling $P_n$. We obtain $T_{\text{MH}} = 0.004821$ which is larger than the estimated quantile $0.002918$; also, we can estimate the p-value of $T_{\text{MH}}$, which is $0.005 \ll \alpha$, suggesting that we reject $H_0 : P_{\text{MH}} =  P_{\text{Null}}$. We apply the same procedure to $P_{\text{Mu}}$ with $n := |P_{\text{Mu}}| =$ 5,731 observations. We obtain $T_{\text{Mu}} = 0.002334$, which is slightly below the estimated quantile $0.002412$, and the p-value of $T_{\text{Mu}}$ is $0.066 > \alpha$, implying that we cannot $H_0 : P_{\text{Mu}} =  P_{\text{Null}}$ at level $\alpha$. Therefore, our testing procedure provides a rigorous framework to detect weak signals that are otherwise indistinguishable, as visualized in Figure \ref{fig:p-values}.

\section{Discussion}
\paragraph{Comparison with existing methods}
Notice that the problem \eqref{eq:displacement} can be viewed as an instance of general testing 
\begin{equation*}
    H_0 : X_1, \ldots, X_n \overset{\mathrm{i.i.d.}}{\sim} P
    \quad \text{vs.} \quad 
    H_1 : X_1, \ldots, X_n \overset{\mathrm{i.i.d.}}{\sim} Q_n,
\end{equation*}
where $Q_n \to P$ as $n \to \infty$; our problem is the case where $Q_n$ is given as displacement interpolation between $P$ and $Q$. Existing approaches to such a general problem characterize detection by utilizing the following (related) notions: likelihood ratio $\dv{Q_n}{P}$ (Section 13.10.1 of \cite{van1996weak}) or Hellinger distance $H^2(P, Q_n)$ (Section 13.1 of \cite{lehmann_2005}). The former approach---often referred to as Le Cam's third lemma---derives limit distributions based on asymptotic normality; the latter characterizes the detection boundary depending on the rate of $n H^2(P, Q_n)$. The first principle behind these methods is clear: quantify a distance/discrepancy between $P, Q_n$ and characterize the limit. Our results also utilize such a first principle: we use weighted Wasserstein distances. Though our method and the existing methods share a similar high-level idea, the existing methods are unsuitable for our problem for several reasons. First, though Le Cam's third lemma applies to any abstract setting under certain conditions, it does not lead to the exact characterization of testing errors unless there is a suitable parametric assumption. Second, the Hellinger distance is not preferred in analyzing the case where $Q_n$ is given as displacement interpolation as its relationship with the interpolation parameter $\epsilon_n$ is not transparent. Moreover, the Hellinger distance-based characterization generally does not calculate testing errors at the detection boundary. On the other hand, our method motivated by weighted Wasserstein distances not only interplays well with displacement interpolation but also provides the exact characterization of testing errors, including the boundary case; moreover, unlike the likelihood ratio test, which requires information on both $P$ and $Q_n$, our test is implementable as long as $P$ is known. 

\paragraph{Lifting technical assumptions}
As remarked in Remark \ref{rmk:lipschitzness}, the Lipschitzness assumption on $G^{-1} \circ F$ is used in the detectable case ($n^{1 / 2} \epsilon_n \to \infty$), but not in the other two cases. Removing the Lipschitzness assumption, there is no restriction on $Q$ as long as its quantile function $G^{-1}$ is bounded. Accordingly, our main results hold even when $Q$ is discrete; in particular, there are cases where our method applies, but Le Cam's third lemma cannot because contiguity is not satisfied. As a concrete example, suppose $P = \mathrm{Unif}[0, 1]$ and $Q = \frac{1}{2} \delta_{0} + \frac{1}{2} \delta_1$. Then, $Q_n := ((1 - \epsilon_n) \mathrm{Id} + \epsilon_n G^{-1} \circ F)_\# P$ is a uniform measure supported on the union of two disjoint intervals $U_n := [0, (1 - \epsilon_n) / 2] \cup [(1 + \epsilon_n) / 2, 1]$. One can verify that $Q_n^{\otimes n}$---the $n$ tensor product of $Q_n$---is not contiguous with respect to $P^{\otimes n}$ because $Q_n^{\otimes n}(U_n^n) = 1$ while $P^{\otimes n}(U_n^n) = (1-\epsilon_n)^n \rightarrow 0$ at the boundary $\epsilon_n \asymp n^{-1 / 2}$. Hence, Le Cam's third lemma cannot be applied to derive the limit distribution for such a case, whereas our method can. Lastly, we briefly discuss Assumption \ref{a:bounded_away}. Though it is a reasonable assumption to impose from a practical viewpoint as discussed in Remark \ref{rmk:on_assumption}, it is natural to ask---from a technical viewpoint---whether such an assumption can be relaxed. To circumvent Assumption \ref{a:bounded_away}, we can use an alternative pivotal test based on the weighted Wasserstein between $\mathrm{Unif}[0, 1]$ and the empirical measure based on $F(X_1), \ldots, F(X_n)$; then, we can show that similar theory holds.

\bibliography{ref.bib}
\bibliographystyle{plainnat}

\appendix
\section{Proofs}
\label{sec:proofs}
\subsection{Proof of Lemma \ref{lem:integrated_quantile_convergence}}
\begin{proof}
    Recall from Lemma \ref{lem:quantile_limit} that $\mathbf{H}_n \leadsto \mathbf{H}$ in $\ell^\infty(0, 1)$. Now, define
    \begin{equation*}
        \ell_m^\infty(0, 1) = \{h \in \ell^\infty(0, 1) : h ~ \text{is measurable}\},
    \end{equation*}
    then one can verify that $\ell_m^\infty(0, 1)$ is a closed subset of the Banach space $(\ell^\infty(0, 1), \|\cdot\|_\infty)$. Also, as the sample path of $\mathbf{H}_n$ is always monotone and the sample path of $\mathbf{H}$ is always continuous, $\mathbf{H}_n$ and $\mathbf{H}$ take values in $\ell_m^\infty(0, 1)$. Now, define a map $\cI \colon \ell_m^\infty(0, 1) \to \R$ by 
    \begin{equation*}
        \cI(h) = \int_{0}^{1} |h(u)|^2 \dd{\omega(u)}.
    \end{equation*}
    We claim that $\cI$ is continuous. To this end, consider a sequence $(h_n)_{n \in \N}$ in $\ell_m^\infty(0, 1)$ and $h \in \ell_m^\infty(0, 1)$ such that $\|h_n - h\|_\infty \to 0$. Then, $\sup_{n \in \N} \|h_n\|_\infty \le M$ for some $M > 0$, hence the dominated convergence theorem shows that $\cI(h_n) \to \cI(h)$. Now, applying Theorem \ref{thm:extended_continuous_mapping} (stated below), we conclude that $\cI(\mathbf{H}_n) \leadsto \cI(\mathbf{H})$, which proves \eqref{eq:convergence1}. Similarly, for each $n \in \N$, define a map $\cE_n \colon \ell_m^\infty(0, 1) \to \R$ by 
    \begin{equation*}
        \cE_n(h) = \int_{0}^{1} |(G^{-1} \circ F)(F^{-1}(u) + n^{-1 / 2} h(u)) - F^{-1}(u)|^2 \dd{\omega(u)}.
    \end{equation*}
    For any sequence $(h_n)_{n \in \N}$ in $\ell_m^\infty(0, 1)$ and $h \in \ell_m^\infty(0, 1)$ such that $\|h_n - h\|_\infty \to 0$, we claim 
    \begin{equation}
        \label{eq:convergence2_pre}
        \cE_n(h_n) \to \int_{0}^{1} |G^{-1}(u) - F^{-1}(u)|^2 \dd{\omega(u)}.
    \end{equation}
    As $G^{-1}$ is continuous almost everywhere on $(0, 1)$, we can see that $G^{-1} \circ F \circ (F^{-1} + n^{-1 / 2} h_n)$ converges to $G^{-1}$ almost everywhere; as $\omega$ is absolutely continuous with respect to the Lebesgue measure, this convergence holds $\omega$-almost everywhere as well. As $G^{-1}$ and $F^{-1}$ are bounded on $(0, 1)$, the dominated convergence theorem shows \eqref{eq:convergence2_pre}. By Theorem \ref{thm:extended_continuous_mapping}, 
    \begin{equation*}
        \cE_n(\mathbf{H}_n) \leadsto \int_{0}^{1} |G^{-1}(u) - F^{-1}(u)|^2 \dd{\omega(u)},
    \end{equation*}
    showing \eqref{eq:convergence2}. Lastly, to prove \eqref{eq:convergence3}, for each $n \in \N$, define a map $\cJ_n \colon \ell_m^\infty(0, 1) \to \R$ by 
    \begin{equation*}
        \cJ_n(h) = \int_{0}^{1} h(u) \cdot \left((G^{-1} \circ F)(F^{-1}(u) + n^{-1 / 2} h(u)) - F^{-1}(u)\right) \dd{\omega(u)}.
    \end{equation*}
    Also, define $\cJ \colon \ell_m^\infty(0, 1) \to \R$ by 
    \begin{equation*}
        \cJ(h) = \int_{0}^{1} h(u) \cdot \left(G^{-1}(u) - F^{-1}(u)\right) \dd{\omega(u)}.
    \end{equation*}
    Similarly to $\cI$, one can verify that $\cJ$ is continuous. Also, for any sequence $(h_n)_{n \in \N}$ in $\ell_m^\infty(0, 1)$ and $h \in \ell_m^\infty(0, 1)$ such that $\|h_n - h\|_\infty \to 0$, we can show that $\cJ_n(h_n) \to \cJ(h)$ by means of the dominated convergence theorem. Therefore, $\cJ_n(\mathbf{H}_n) \leadsto \cJ(\mathbf{H})$ holds by Theorem \ref{thm:extended_continuous_mapping}, hence \eqref{eq:convergence3} holds.
\end{proof}

\begin{theorem}[Extended Continuous Mapping Theorem]
    \label{thm:extended_continuous_mapping}
    Let $\cH$ and $\cK$ be metric spaces. Let $\cH_0$ be a Borel subset of $\cH$ and consider a measurable map $\cI \colon \cH_0 \to \cK$. Suppose there is a sequence $(\cI_n)_{n \in \N}$ of maps from $\cH_0$ to $\cK$ satisfying the following: for any sequence $(h_n)_{n \in \N}$ in $\cH_0$ converging to some $h \in \cH_0$, the sequence $(\cI_n(h_n))_{n \in \N}$ converges to $\cI(h)$ in $\cK$. If a sequence $(H_n)_{n \in \N}$ of random elements in $\cH$ converges weakly to a tight measurable random element $H$ in $\cH$, where both $H_n$ and $H$ take values in $\cH_0$, the sequence $(\cI_n(H_n))_{n \in \N}$ of random elements in $\cK$ converges weakly to the random element $\cI(H)$ in $\cK$.
\end{theorem}

\begin{remark}
    \textup{Theorem \ref{thm:extended_continuous_mapping} is adapted from Theorem 1.11.1 of \cite{van1996weak}, where the latter uses separability instead of tightness. As tightness implies separability, we have stated Theorem \ref{thm:extended_continuous_mapping} with tightness, which is sufficient in our analysis.}
\end{remark}

\subsection{Proof of Theorem \ref{thm:asymptotics}}
\begin{proof}
    For each $n \in \N$, let $\phi_n = (1 - \epsilon_n) \mathrm{Id} + \epsilon_n G^{-1} \circ F$. As we are concerned with the limit distribution of $n W_{2, \omega}^2(P_n, P)$, we may assume $(X_n)_{n \in \N}$ is i.i.d.\ from $P$ and let $P_n$ be the empirical measure based on $\phi_n(X_1), \ldots, \phi_n(X_n)$ as $\phi_n(X_1), \ldots, \phi_n(X_n)$ also follow $H_1^{(n)}$. Now, let $F_n^{-1}$ be the empirical quantile function based on $X_1, \ldots, X_n$, then 
    \begin{equation*}
        \begin{split}
            \phi_n \circ F_n^{-1} - F^{-1} & = (1 - \epsilon_n) F_n^{-1} + \epsilon_n G^{-1} \circ F \circ F_n^{-1} - F^{-1} \\
            & = (1 - \epsilon_n) n^{-1 / 2} \mathbf{H}_n + \epsilon_n (G^{-1} \circ F \circ (F^{-1} + n^{-1 / 2} \mathbf{H}_n) - F^{-1}), 
        \end{split}
    \end{equation*}
    where $\mathbf{H}_n := \sqrt{n} (F_n^{-1} - F^{-1})$. Hence, using $A_n, B_n, C_n$ defined in \eqref{eq:A_n}, \eqref{eq:B_n}, \eqref{eq:C_n}, respectively, we have
    \begin{equation*}
        W_{2, \omega}^2(P_n, P) = (1 - \epsilon_n)^2 n^{-1} A_n + \epsilon_n^2 B_n + 2 (1 - \epsilon_n) n^{-1 / 2} \epsilon_n C_n.
    \end{equation*} 
    \textbf{Case I}: Suppose $n^{1 / 2} \epsilon_n \to 0$. Recall that we have shown in Lemma \ref{lem:integrated_quantile_convergence} that $A_n$, $B_n$, and $C_n$ are weakly convergent. Note that
    \begin{equation*}
        n W_{2, \omega}^2(P_n, P) = (1 - \epsilon_n)^2 A_n + (n^{1 / 2} \epsilon_n)^2 B_n + 2 (1 - \epsilon_n) (n^{1 / 2} \epsilon_n) C_n.
    \end{equation*}
    As $(n^{1 / 2} \epsilon_n)^2 B_n, (n^{1 / 2} \epsilon_n) C_n \leadsto 0$, by Slutsky's theorem,
    \begin{equation*}
        n W_{2, \omega}^2(P_n, P) \leadsto \lim_{n \to \infty} A_n = \int_{0}^{1} |\mathbf{H}(u)|^2 \dd{\omega(u)}.
    \end{equation*}
    \textbf{Case II}: Suppose $n^{1 / 2} \epsilon_n \to \infty$ and note that
    \begin{equation*}
        \begin{split}
            n^{1 / 2} \epsilon_n^{-1} \left(W_{2, \omega}^2(P_n, P) - \epsilon_n^2 W_{2, \omega}^2(P, Q)\right) & = (1 - \epsilon_n)^2 (n^{1 / 2} \epsilon_n)^{-1} A_n \\
            & \quad + n^{1 / 2} \epsilon_n (B_n - W_{2, \omega}^2(P, Q)) + 2 (1 - \epsilon_n) C_n.
        \end{split}
    \end{equation*}
    First, notice that $(n^{1 / 2} \epsilon_n)^{-1} A_n \leadsto 0$. We claim $n^{1 / 2} \epsilon_n (B_n - W_{2, \omega}^2(P, Q)) \leadsto 0$. To this end, observe that
    \begin{equation*}
        \begin{split}
            & B_n - W_{2, \omega}^2(P, Q) \\
            & = \int_{0}^{1} |(G^{-1} \circ F)(F^{-1}(u) + n^{-1 / 2} \mathbf{H}_n(u)) - F^{-1}(u)|^2 \dd{\omega(u)} \\
            & \quad - \int_{0}^{1} |G^{-1}(u) - F^{-1}(u)|^2 \dd{\omega(u)} \\
            & = \int_{0}^{1} |(G^{-1} \circ F)(F^{-1}(u) + n^{-1 / 2} \mathbf{H}_n(u)) - G^{-1}(u)|^2 \dd{\omega(u)} \\
            & \quad + 2 \int_{0}^{1} \left((G^{-1} \circ F)(F^{-1}(u) + n^{-1 / 2} \mathbf{H}_n(u)) - G^{-1}(u)\right) \cdot \left(G^{-1}(u) - F^{-1}(u)\right) \dd{\omega(u)} \\
            & =: B_n^o + 2 B_n^\ast.
        \end{split}
    \end{equation*}
    Let $L$ be the Lipschitz constant of $G^{-1} \circ F$, then
    \begin{equation*}
        n^{1 / 2} \epsilon_n |B_n^\ast| 
        \le n^{1 / 2} \epsilon_n \int_{0}^{1} L n^{-1 / 2} |\mathbf{H}_n(u)| \cdot |G^{-1}(u) - F^{-1}(u)| \dd{\omega(u)}
        \le L \epsilon_n \sqrt{A_n} W_{2, \omega}^2(P, Q),
    \end{equation*}
    where $\epsilon_n \sqrt{A_n} \leadsto 0$ by Lemma \ref{lem:integrated_quantile_convergence}, hence $n^{1 / 2} \epsilon_n B_n^\ast \leadsto 0$. Similarly, 
    \begin{equation*}
        n^{1 / 2} \epsilon_n |B_n^0| \le n^{1 / 2} \epsilon_n \int_{0}^{1} L^2 n^{-1} |\mathbf{H}_n(u)|^2 \dd{\omega(u)} = L^2 n^{-1 / 2} \epsilon_n A_n \leadsto 0,
    \end{equation*}
    hence $n^{1 / 2} \epsilon_n B_n^o \leadsto 0$. Therefore, applying Slutsky's theorem, we have
    \begin{equation*}
        \begin{split}
            n^{1 / 2} \epsilon_n^{-1} \left(W_{2, \omega}^2(P_n, P) - \epsilon_n^2 W_{2, \omega}^2(P, Q)\right) & \leadsto \lim_{n \to \infty} 2 (1 - \epsilon_n) C_n \\
            & = 2 \int_{0}^{1} \mathbf{H}(u) \cdot \left(G^{-1}(u) - F^{-1}(u)\right) \dd{\omega(u)}.            
        \end{split}
    \end{equation*}
    \textbf{Case III}: Suppose $n^{1 / 2} \epsilon_n = \gamma$, then
    \begin{equation*}
        n W_{2, \omega}^2(P_n, P) = (1 - \epsilon_n)^2 A_n + \gamma^2 B_n + 2 (1 - \epsilon_n) \gamma C_n = (1 - \epsilon_n) \cdot \left((1 - \epsilon_n) A_n + 2 \gamma C_n\right) + \gamma^2 B_n.
    \end{equation*}
    We claim that 
    \begin{equation*}
        (1 - \epsilon_n) A_n + 2 \gamma C_n \leadsto \int_{0}^{1} |\mathbf{H}(u)|^2 \dd{\omega(u)} + 2 \gamma \int_{0}^{1} \mathbf{H}(u) \cdot \left(G^{-1}(u) - F^{-1}(u)\right) \dd{\omega(u)}.           
    \end{equation*}
    We apply the same argument used in the proof of Lemma \ref{lem:integrated_quantile_convergence}; for each $n \in \N$, define $\cK_n \colon \ell_m^\infty(0, 1) \to \R$ by 
    \begin{equation*}
        \begin{split}
            \cK_n(h) & = (1 - \epsilon_n) \int_{0}^{1} |h(u)|^2 \dd{\omega(u)} \\
            & \quad + 2 \gamma \int_{0}^{1} h(u) \cdot \left((G^{-1} \circ F)(F^{-1}(u) + n^{-1 / 2} h(u)) - F^{-1}(u)\right) \dd{\omega(u)}.            
        \end{split}
    \end{equation*}
    In the proof of Lemma \ref{lem:integrated_quantile_convergence}, we have shown that 
    \begin{equation*}
        \cK_n(h_n) \to \int_{0}^{1} |h(u)|^2 \dd{\omega(u)} + 2 \gamma \int_{0}^{1} h(u) \cdot \left(G^{-1}(u) - F^{-1}(u)\right) \dd{\omega(u)}
    \end{equation*}
    for any sequence $(h_n)_{n \in \N}$ in $\ell_m^\infty(0, 1)$ and $h \in \ell_m^\infty(0, 1)$ such that $\|h_n - h\|_\infty \to 0$. Therefore, the extended continuous mapping theorem shows that 
    \begin{equation*}
        (1 - \epsilon_n) A_n + 2 \gamma C_n = \cK_n(\mathbf{H}_n) \leadsto \int_{0}^{1} |\mathbf{H}(u)|^2 \dd{\omega(u)} + 2 \gamma \int_{0}^{1} \mathbf{H}(u) \cdot \left(G^{-1}(u) - F^{-1}(u)\right) \dd{\omega(u)}.     
    \end{equation*}
    Hence, by Slutsky's theorem, 
    \begin{equation*}
        \begin{split}
            & n W_{2, \omega}^2(P_n, P) \\
            & = (1 - \epsilon_n) \cdot \left((1 - \epsilon_n) A_n + 2 \gamma C_n\right) + \gamma^2 B_n \\
            & \leadsto \int_{0}^{1} |\mathbf{H}(u)|^2 \dd{\omega(u)} + 2 \gamma \int_{0}^{1} \mathbf{H}(u) \cdot \left(G^{-1}(u) - F^{-1}(u)\right) \dd{\omega(u)} + \gamma^2 W_{2, \omega}^2(P, Q).
        \end{split}
    \end{equation*}
\end{proof}

\subsection{Proof of Theorem \ref{thm:asymptotic_errors}}
\begin{proof}
    For the asymptotic Type I error, we invoke Proposition \ref{prop:limit_under_H0}: letting $P_n$ be the empirical measure based on $X_1, \ldots, X_n$ from $H_0^{(n)}$,
    \begin{equation*}
        \text{the asymptotic Type I error} = \lim_{n \to \infty} \mathbb{P}(n W_{2, \omega}^2(P_n, P) > C_\alpha) = \mathbb{P}(A > C_\alpha) = 1 - \Psi(C_\alpha) = \alpha,
    \end{equation*}
    where $A$ is the limit distribution of $n W_{2, \omega}^2(P_n, P)$ under $H_0^{(n)}$, namely, $A =$ \eqref{eq:limit_distribution_null}. Here, weak convergence implies the above limit as the open set $(C_\alpha, \infty)$ is a continuity set of $\Psi$. Therefore, the asymptotic Type I error is $\alpha$. 

    Next, we compute the asymptotic Type II error: assuming $P_n$ be the empirical measure based on $X_1, \ldots, X_n$ from $H_1^{(n)}$,
    \begin{equation*}
        \text{the asymptotic Type II error} = \lim_{n \to \infty} \mathbb{P}(n W_{2, \omega}^2(P_n, P) \le C_\alpha).
    \end{equation*}
    If $n^{1 / 2} \epsilon_n \to 0$, we have shown in Theorem \ref{thm:asymptotics} that the limit distribution of $n W_{2, \omega}^2(P_n, P)$ is exactly $A$, hence
    \begin{equation*}
        \lim_{n \to \infty} \mathbb{P}(n W_{2, \omega}^2(P_n, P) \le C_\alpha) = \mathbb{P}(A \le C_\alpha) = \Psi(C_\alpha) = 1 - \alpha.
    \end{equation*}
    If $n^{1 / 2} \epsilon_n \to \infty$, we apply Slutsky's theorem by multiplying $n^{-1 / 2} \epsilon_n^{-1}$ to both sides of \eqref{eq:beta<1/2}, which yields
    \begin{equation*}
        \epsilon_n^{-2} \left(W_{2, \omega}^2(P_n, P) - \epsilon_n^2 W_{2, \omega}^2(P, Q)\right) \leadsto 0.
    \end{equation*}
    In other words, $\epsilon_n^{-2} W_{2, \omega}^2(P_n, P)$ converges to a constant $W_{2, \omega}^2(P, Q) > 0$ under $H_1^{(n)}$, hence 
    \begin{equation*}
        \lim_{n \to \infty} \mathbb{P}(n W_{2, \omega}^2(P_n, P) \le C_\alpha) = \lim_{n \to \infty} \mathbb{P}(\epsilon_n^{-2} W_{2, \omega}^2(P_n, P) \le n^{-1} \epsilon_n^{-2} C_\alpha) = 0.
    \end{equation*}
    If $n^{1 / 2} \epsilon_n = \gamma > 0$, the limit distribution of $T_n := n W_{2, \omega}^2(P_n, P) - \gamma^2 W_{2, \omega}^2(P, Q)$ under $H_1^{(n)}$ is \eqref{eq:limit_distribution_alt} by \eqref{eq:beta=1/2}, hence 
    \begin{equation*}
        \lim_{n \to \infty} \mathbb{P}(n W_{2, \omega}^2(P_n, P) \le C_\alpha) = \lim_{n \to \infty} \mathbb{P}(T_n \le C_\alpha - \gamma^2 W_{2, \omega}^2(P, Q)) = \Psi_\gamma(C_\alpha - \gamma^2 W_{2, \omega}^2(P, Q)).
    \end{equation*}
	Here, weak convergence implies the above limit as the set $(-\infty, C_\alpha]$ is a continuity set of $\Psi_\gamma$.
\end{proof}
\begin{remark}
    \textup{As noted in Remark \ref{rmk:main}, we may replace the detection boundary $n^{1 / 2} \epsilon_n = \gamma > 0$ with $\lim_{n \to \infty} n^{1 / 2} \epsilon_n = \gamma > 0$ by modifying the proofs of Theorem \ref{thm:asymptotics} and \ref{thm:asymptotic_errors}.}
\end{remark}

\end{document}